\documentclass[12pt]{article}
\usepackage[margin=1.1in]{geometry}

\usepackage{graphicx}
\usepackage{amsmath,amsthm,amssymb}
\usepackage{amsthm,amssymb}

\newtheorem{thm}{Theorem}[section]

\newtheorem{rem}[thm]{Remark}

\usepackage{hyperref}
\hypersetup{colorlinks=true}

\numberwithin{equation}{section}

\newcommand{\reviewA}[1]{\textcolor[rgb]{0.00,0.00,0.00}{#1}}

\newcommand{\reviewC}[1]{\textcolor[rgb]{0.00,0.00,0.00}{#1}}

\newcommand{\rmd}{\mathrm{d}}
\newcommand{\rme}{\mathrm{e}}
\newcommand{\rmi}{\mathbf{i}}

\newcommand{\eref}[1]{\eqref{#1}}
\newcommand{\Sref}[1]{Section \ref{#1}}
\newcommand{\Fref}[1]{Figure \ref{#1}}
\newcommand{\eqalign}[1]{\begin{aligned}#1\end{aligned}}

\begin{document}

\title{Increasing stability in the linearized inverse Schr\"{o}dinger potential problem with power type nonlinearities
\footnote{ %
S. Lu is supported by NSFC (No.11925104), Science and Technology Commission of Shanghai Municipality (19XD1420500, 21JC1400500). %
M. Salo is supported by the Academy of Finland (Finnish Centre of Excellence in Inverse Modelling and Imaging, grant 284715) and by the European Research Council under Horizon 2020 (ERC CoG 770924). %
B. Xu is supported by NSFC (No.12171301 and No.11801351). %
}}




\author{%
Shuai Lu%
\footnote{Shanghai Key Laboratory for Contemporary Applied Mathematics, Key Laboratory of Mathematics for Nonlinear Sciences and School of Mathematical Sciences, Fudan University, Shanghai, China. Email: \textsf{slu@fudan.edu.cn}}, %
Mikko Salo%
\footnote{Department of Mathematics and Statistics, University of Jyv\"{a}skyl\"{a}, Jyv\"{a}skyl\"{a}, Finland. Email: \textsf{mikko.j.salo@jyu.fi}}, %
Boxi Xu%
\footnote{Author to whom any correspondence should be addressed. School of Mathematics, Shanghai University of Finance and Economics, Shanghai, China. Email: \textsf{xu.boxi@mail.sufe.edu.cn}}%
}

\maketitle

\begin{abstract}
We consider increasing stability in the inverse Schr\"{o}dinger potential problem with power type nonlinearities at a large wavenumber. Two linearization approaches, with respect to small boundary data and small potential function, are proposed and their performance on the inverse Schr\"{o}dinger potential problem is investigated. It can be observed that higher order linearization for small boundary data can provide an increasing stability for an arbitrary power type nonlinearity term if the wavenumber is chosen large. Meanwhile, linearization with respect to the potential function leads to increasing stability for a quadratic nonlinearity term, which highlights the advantage of nonlinearity in solving the inverse Schr\"{o}dinger potential problem. Noticing that both linearization approaches can be numerically approximated, we provide several reconstruction algorithms for the quadratic and general power type nonlinearity terms, where one of these algorithms is designed based on boundary measurements of multiple wavenumbers. Several numerical examples shed light on the efficiency of our proposed algorithms.
\end{abstract}





\section{Introduction}
\subsection{Background}

The inverse Schr\"{o}dinger potential problem arises from electrical impedance tomography (EIT) \cite{C1980} and has attracted much attention both theoretically and computationally. In a general setting, we can formulate the following Schr\"{o}dinger equation
\begin{equation}\label{eq_generalInverSchrodinger}
\left\{~
\eqalign{
\Delta u + k^{2} u - c(x) u = 0 &\quad \textrm{in\ } \Omega \subset \mathbb{R}^{n}, \\
u = g_{0} &\quad \textrm{on\ } \partial\Omega,
}
\right.
\end{equation}
where, throughout the article, $\Omega \subset \mathbb{R}^{n}$ is assumed to be a bounded open domain with smooth boundary $\partial\Omega$ and dimension $n \geq 2$. The inverse Schr\"{o}dinger potential problem is to identify the unknown potential function $c(x)$ from many boundary measurements or the Dirichlet-to-Neumann map defined below. A classical result in \cite{A1988} shows that if the wavenumber $k = 0$ in \eref{eq_generalInverSchrodinger} the stability of the inverse Schr\"{o}dinger potential problem is logarithmic. When the wavenumber is sufficiently large, increasing stability with respect to the wavenumber $k$ has been observed and well documented, starting with \cite{I2011} and with many further results given in \cite{IW2014, ILW2016, ILX2020} for \eref{eq_generalInverSchrodinger} or its linearized form. These results are often stated as stability estimates involving a H\"{o}lder term and a logarithmic term which goes to zero as the wavenumber goes to infinity. An alternative way to observe increasing stability is to note that one can determine the Fourier transform of the unknown coefficient in a stable way for a range of frequencies, and that this range increases with the wavenumber. We note that these increasing stability results have also been verified both theoretically and numerically in other inverse source, obstacle or medium problems where we refer to \cite{BLT2010, BT2010, NUW2013, BG2015, BLLT2015, BLRX2015, CIL2016, ZZ2017, IL2018, KKK2018, BLZ2020, BT2020} and references therein.

There have also been several recent works on inverse problems for nonlinear elliptic equations. In such problems, it has been observed that higher order linearizations of the nonlinear Dirichlet-to-Neumann map carry information about the unknown coefficients. This method allows one to exploit nonlinear effects in order to obtain better results than those that are currently known for corresponding linear equations. The higher order linearization method goes back to \cite{KLU2018} in the hyperbolic case and to \cite{FO2020, LLLS_I} in the elliptic case. The method has been further applied to more general equations and partial data problems. See \cite{KU2019, KKU2020, LLST2020, CFKKU2021, LLLS_II} for a selection of recent results.

This article studies possible improvements in stability properties of inverse problems for nonlinear Schr\"{o}dinger type equations with a large wavenumber.
More specifically, we study the inverse Schr\"{o}dinger potential problem with an arbitrary power type nonlinearity term and discuss its unique determination, increasing stability and numerical reconstruction algorithms. In particular, we consider the problem of recovering the potential function $c(x)$, defined in $\Omega \subset \mathbb{R}^{n}$, in the following nonlinear Schr\"{o}dinger equation, with an integer $m \geq 2$ denoting the nonlinearity index,
\begin{equation}\label{eq_HelmholtzEqmain}
\left\{~
\eqalign{
\Delta u + k^{2} u - c(x) u^{m} = 0 &\quad \textrm{in\ } \Omega, \\
u = g_{0} &\quad \textrm{on\ } \partial\Omega,
}
\right.
\end{equation}
from many boundary measurements. Here, we assume that the squared wavenumber $k^{2}$ is sufficiently large, $0$ is not a Dirichlet eigenvalue of $\Delta + k^{2}$ in $\Omega$ \reviewC{and the Dirichlet boundary data $g_{0}$ is sufficiently small}. Meanwhile, by assuming that $c(x)$ is compactly supported in $\Omega$, the well-posedness of the forward problem \eref{eq_HelmholtzEqmain} can be verified following the variational framework developed in \cite[Theorem 1]{EW2014}. Thus, the boundary measurements can be given by the nonlinear Dirichlet-to-Neumann (DtN) map
\begin{equation}\label{eq_DtN}
\Lambda_{c}: g_{0} \mapsto \partial_{\nu} u \quad \textrm{on\ } \partial\Omega.
\end{equation}
The precise definition of $\Lambda_{c}$ and its two linearized forms $D^{m}_{0} \Lambda_{c}$, $\Lambda'_{c}$ will be specified later.
\reviewC{When the wavenumber $k = 0$ in (\ref{eq_HelmholtzEqmain}), unique identification of the potential function $c(x)$ has been provided in \cite{LLLS_II} by measurement of the DtN map in (\ref{eq_DtN}) and its linearized form $D^{m}_{0} \Lambda_{c}$. In current work, we particularly focus on the increasing stability estimate for (\ref{eq_HelmholtzEqmain}) under two different linearized forms of $\Lambda_{c}$.}


\subsection{Linearization approaches}\label{se1.2}

To solve the nonlinear inverse Schr\"{o}dinger potential problem stably, we implement linearization approaches and discuss recovery of the potential function by the linearized DtN maps accordingly. In this subsection, we briefly overview two linearization approaches with respect to small boundary data and small potential function, \reviewC{which have been studied in linear and nonlinear elliptic inverse problems, for instance in \cite{C1980, LLLS_I}, when $k = 0$ in \eref{eq_generalInverSchrodinger} and \eref{eq_HelmholtzEqmain}.}

To treat elliptic equations with power type nonlinearities, a novel linearization approach with respect to small boundary data has recently been discussed in \cite{FO2020, LLLS_I} \reviewA{and been extended to a fractional nonlinearity index in \cite{LLST2020}}.
We briefly introduce its extension to the nonlinear Schr\"{o}dinger potential problem \eref{eq_HelmholtzEqmain} below. Assume that $c \in C^{\alpha}(\overline{\Omega})$ for some $\alpha$ with $0 < \alpha < 1$, and $0$ is not a Dirichlet eigenvalue of $\Delta + k^{2}$ in $\Omega$. By \cite[Proposition 2.1]{LLST2020}, we can find a constant $\tau > 0$ such that for any Dirichlet boundary value $f$ in $U_{\tau} := \{ f \in C^{2,\alpha}(\partial\Omega) \,:\, \|f\|_{C^{2,\alpha}(\partial\Omega)} \leq \tau \}$, there is a unique small solution $u \in C^{2,\alpha}(\overline{\Omega})$ and $u|_{\partial\Omega} = f$. The nonlinear DtN map in the H\"{o}lder spaces is defined by
\begin{equation*}
\Lambda_{c}: U_{\tau} \subset C^{2,\alpha}(\partial\Omega) \to C^{1,\alpha}(\partial\Omega), \qquad
f \mapsto \partial_{\nu} u |_{\partial\Omega}.
\end{equation*}%
Let $\varepsilon = (\varepsilon_{1}, \ldots, \varepsilon_{m})$ where each $\varepsilon_{j} > 0$ is small, and consider the solution $u_{\varepsilon}$ corresponding to the Dirichlet boundary value
\begin{equation*}
f_{\varepsilon} = \varepsilon_{1} f_{1} + \ldots + \varepsilon_{m} f_{m}.
\end{equation*}
By \cite[Proposition 2.1]{LLST2020} the solution $u_{\varepsilon}$ depends smoothly on the parameters $\varepsilon_{j}$. We may thus differentiate the equation
\begin{equation}\label{eq_Linearization1}
\Delta u_{\varepsilon} + k^{2} u_{\varepsilon} - c(x) u_{\varepsilon}^{m} = 0 \quad \textrm{in\ } \Omega, \qquad
u_{\varepsilon} |_{\partial\Omega} = f_{\varepsilon}
\end{equation}
with respect to the parameters $\varepsilon_{j}$. Writing $v_{j} = \partial_{\varepsilon_{j}} u_{\varepsilon} |_{\varepsilon = 0}$, we observe that $v_{j}$ is the unique solution of
\begin{equation*}
\Delta v_{j} + k^{2} v_{j} = 0 \quad \textrm{in\ } \Omega, \qquad
v_{j} |_{\partial\Omega} = f_{j}.
\end{equation*}
Similarly, applying $\partial_{\varepsilon_{1}} \cdots \partial_{\varepsilon_{m}}$ to the equation \eref{eq_Linearization1} and setting $\varepsilon = 0$, we can define $w = \partial_{\varepsilon_{1}} \cdots \partial_{\varepsilon_{m}} u_{\varepsilon} |_{\varepsilon = 0}$ which solves the equation
\begin{equation}\label{mth_equation_difference}
\Delta w + k^{2} w = (m!) c(x) v_{1} \cdots v_{m} \quad \textrm{in\ } \Omega, \qquad
w |_{\partial\Omega} = 0.
\end{equation}
Moreover, the Neumann boundary data can be obtained in form of
\begin{equation}\label{eq_HighOrderLinearDtN}
\eqalign{
\partial_{\nu} w |_{\partial\Omega}
&= \partial_{\varepsilon_{1}} \cdots \partial_{\varepsilon_{m}} ( \partial_{\nu} u_{\varepsilon} |_{\partial\Omega} ) |_{\varepsilon = 0}
= \partial_{\varepsilon_{1}} \cdots \partial_{\varepsilon_{m}} \Lambda_{c} (f_{\varepsilon}) |_{\varepsilon = 0} \\
&= D^{m}_{0} \Lambda_{c} (f_{1}, \ldots, f_{m})
}
\end{equation}
where $D^{m}_{0}$ denotes the $m$th Fr\'{e}chet derivative at $0$ considered as an $m$-linear form. If we integrate the equation \eref{mth_equation_difference} against another function $v_{m+1}$ solving
\begin{equation*}
\Delta v_{m+1} + k^{2} v_{m+1} = 0 \quad \textrm{in\ } \Omega, \qquad
v_{m+1} |_{\partial\Omega} = f_{m+1},
\end{equation*}
we obtain a Calder\'{o}n or Alessandrini type identity
\begin{equation}\label{eq_nonlinearsmallboundary_equality}
(m!) \int_{\Omega} c(x) v_{1} \cdots v_{m} v_{m+1} ~\rmd x
= \int_{\partial\Omega} D^{m}_{0} \Lambda_{c} (f_{1}, \ldots, f_{m}) f_{m+1} ~\rmd S
\end{equation}
which will be revisited later.

Noticing that the $m$th Fr\'{e}chet derivative $D^{m}_{0} \Lambda_{c} (f_{1}, \ldots, f_{m})$ is numerically hard to obtain, we further consider the case when $c(x)$ is small compared to the wavenumber and study the linearization approach with respect to the potential function as investigated in the linear Schr\"{o}dinger potential problem in \cite{ILX2020}.
Taking the asymptotic expansion with respect to the potential function $c(x)$, we have
\begin{equation}\label{eq_asymILX}
u = u_{0} + u_{1} + u_{2} + \ldots
\end{equation}
where the remaining ``$\ldots$" denotes the ``higher" order term and following subproblems are satisfied such that
\begin{equation*}
\eqalign{
& \Delta u_{0} + k^{2} u_{0} = 0, \\
& \Delta u_{1} + k^{2} u_{1} = c(x) u_{0}^{m}, \\
& \Delta u_{2} + k^{2} u_{2} = m c(x) u_{0}^{m-1} u_{1}.
}
\end{equation*}
This shows that $u_{0}$ satisfies the Helmholtz equation $\Delta u_{0} + k^{2} u_{0} = 0$ in $\Omega$ and the first-order expansion term $u_{1}$ satisfies
\begin{equation}\label{eq_quadratic_u1}
\Delta u_{1} + k^{2} u_{1} = c(x) u_{0}^{m} \quad \textrm{in\ } \Omega.
\end{equation}
When $u_{0} |_{\partial\Omega} = g_{0}$ and $u_{1} |_{\partial\Omega} = g_{1} \equiv 0$, the linearized DtN map $\Lambda'_{c}$ is formally defined by
\begin{equation}\label{eq_LinearDtN}
\Lambda'_{c}: g_{0} \mapsto \partial_{\nu} u_{1} \quad \textrm{on\ } \partial\Omega.
\end{equation}
Note that $\Lambda'_{c}$ is actually a nonlinear map, since it corresponds to linearization with respect to the potential.
Multiplying the above equation \eref{eq_quadratic_u1} from both sides with another $\varphi$ solving $\Delta \varphi + k^{2} \varphi = 0$ in $\Omega$, we thus obtain another Calder\'{o}n or Alessandrini type identity
\begin{equation}\label{eq_calderonILX}
\int_{\Omega} c(x) u_{0}^{m} \varphi ~\rmd x
= \int_{\partial\Omega} \partial_{\nu} u_{1} \varphi ~\rmd S
\end{equation}
which will also be revisited later.

In current article, we consider the following problem:

\begin{quote}
{\bf Recover the potential function $c(x)$ from the linearized DtN map $D^{m}_{0} \Lambda_{c}$ or $\Lambda'_{c}$}.
\end{quote}

The first main result shows that from the knowledge of the $m$th Fr\'{e}chet derivative $D^{m}_{0} \Lambda_{c}$, one can determine the Fourier transform $\mathcal{F}[c](\xi)$ of $c$ in a stable way for frequencies $|\xi| \leq (m+1)k$. Thus the range of frequencies that can be determined stably increases both with respect to the wavenumber $k$ and the nonlinearity index $m$. However, determining $D^{m}_{0} \Lambda_{c}$ from $\Lambda_{c}$ becomes numerically very difficult when $m$ increases. The second main result considers the case where the potential function is small compared to the wavenumber. In this case we consider the linearization $\Lambda'_{c}$. We show that in the quadratic case where $m = 2$, from the knowledge of $\Lambda'_{c}$ one can stably determine $\mathcal{F}[c](\xi)$ for frequencies $|\xi| \leq 3k$. This is in contrast with the linear case where one can only determine frequencies $|\xi| \leq 2k$ stably \cite{ILX2020}. Thus in both main results above, the nonlinearity leads to improved stability properties in a certain sense. The theoretical stability results are confirmed by numerical results given in the end of the article.

The article is organized as follows. In \Sref{se2} we show that the linearized DtN map $D^{m}_{0} \Lambda_{c}$ provides a uniform increasing stability where the range of frequencies that can be determined stably increases with respect to $k$ and $m$. On the other hand, $\Lambda'_{c}$ only yields the uniqueness of the potential function $c(x)$ in the general setting $m \geq 2$.
In \Sref{se3} we further explore the linearized DtN map $\Lambda'_{c}$ for the inverse Schr\"{o}dinger potential problem with a quadratic nonlinearity term. By calibrating the identity \eref{eq_calderonILX} carefully, we verify an improved increasing stability for the specific inverse Schr\"{o}dinger potential problem with a quadratic nonlinearity term.
Noticing that both linearized DtN maps $D^{m}_{0} \Lambda_{c}$ and $\Lambda'_{c}$ can be numerically approximated, we extend the reconstruction algorithm in \cite{ILX2020} to the inverse Schr\"{o}dinger potential problem with quadratic and general nonlinearity terms in \Sref{se4}, respectively. We note that one of these reconstruction algorithms is realized by the linearized DtN map $\Lambda'_{c}$ with multiple wavenumbers. In the same \Sref{se4} we provide some numerical examples and extended discussion verifying the efficiency of our proposed algorithms.


\section{Linearized inverse Schr\"{o}dinger potential problem with an arbitrary power type nonlinearity term}\label{se2}

In this section, we investigate the linearized inverse Schr\"{o}dinger potential problem with an arbitrary power type nonlinearity term provided with the linearized DtN map $D^{m}_{0} \Lambda_{c}$ or $\Lambda'_{c}$. The analysis is based on the Calder\'{o}n type identities \eref{eq_nonlinearsmallboundary_equality} and \eref{eq_calderonILX}.

\reviewC{For the map $D^{m}_{0} \Lambda_{c}$ with $k = 0$, \cite{LLLS_I} has verified that by linearizing the small boundary data, the stability estimate for the inverse potential problem is logarithmic, which is consistent with the classical result in EIT \cite{A1988}.} In the current section, we verify that by constructing an appropriate set of complex exponential solutions, there will be improved stability when the wavenumber is large.

Recall the identity \eref{eq_nonlinearsmallboundary_equality},
\begin{equation*}
(m!) \int_{\Omega} c(x) v_{1} \cdots v_{m} v_{m+1} ~\rmd x
= \int_{\partial\Omega} D^{m}_{0} \Lambda_{c} (f_{1}, \ldots, f_{m}) f_{m+1} ~\rmd S.
\end{equation*}
Here $v_{j}$ solve $\Delta v_{j} + k^{2} v_{j} = 0$ in $\Omega$ with $v_{j} |_{\partial\Omega} = f_{j}$. To derive the stability estimate, we rely on the above identity and observe that
\begin{equation}\label{eq_linearization1_estimate}
\left| \int_{\Omega} c(x) v_{1} \cdots v_{m} v_{m+1} ~\rmd x \right|
\leq \frac{\epsilon}{m!} \left( \prod_{j=1}^{m} \| f_{j} \|_{C^{2,\alpha}(\partial\Omega)} \right) \| f_{m+1} \|_{L^{2}(\partial\Omega)}
\end{equation}
where we define $\epsilon := \sup_{\|f_{j}\|_{C^{2,\alpha}(\partial\Omega)} \leq 1} \left\|D^{m}_{0} \Lambda_{c} (f_{1}, \ldots, f_{m})\right\|_{L^{2}(\partial\Omega)}$. Thus \eref{eq_linearization1_estimate} further yields the inequality
\begin{equation}\label{mth_derivative_identity_estimate}
\left| \int_{\Omega} c(x) v_{1} \cdots v_{m+1} ~\rmd x \right|
\leq \frac{\epsilon}{m!} \prod_{j=1}^{m+1} \| v_{j} \|_{C^{2,\alpha}(\overline{\Omega})}.
\end{equation}

Let $\mathcal{F}[c](\xi)$ denote the Fourier transform of $c$ (extended by zero outside $\Omega$) at a frequency $\xi \in \mathbb{R}^{n}$. The following result shows that frequencies $|\xi| \leq (m+1)k$ can be recovered in a Lipschitz stable way from the knowledge of the linearized map $D^{m}_{0} \Lambda_{c}$.

\begin{thm}
Let $m \geq 2$ be an integer, $k \geq 1$, and assume that $|\xi| \leq (m+1)k$. Then
\begin{equation*}
| \mathcal{F}[c](\xi) |
\leq \frac{\epsilon}{m!} \left( 3(1+k^6) \right)^{\frac{m+1}{2}}.
\end{equation*}
\end{thm}

\begin{proof}
We first claim that if $\ell \geq 2$ is an integer, then for any $\eta \in \mathbb{R}^{n}$ with $|\eta| \leq \ell$ there are unit vectors $\omega_{1}, \ldots, \omega_{\ell} \in \mathbb{R}^{n}$ such that
\begin{equation*}
\sum_{j=1}^{\ell} \omega_{j} = \eta.
\end{equation*}
This can be proved by induction. When $\ell = 2$ and $|\eta| \leq 2$, we may choose
\begin{equation*}
\omega_{1} = \frac{\eta}{2} + \sqrt{1-\frac{|\eta|^{2}}{4}}\, \omega, \qquad
\omega_{2} = \frac{\eta}{2} - \sqrt{1-\frac{|\eta|^{2}}{4}}\, \omega,
\end{equation*}
where $\omega$ is any unit vector orthogonal to $\eta$. We make the induction hypothesis that the claim holds for some $\ell \geq 2$. Let $\eta$ be a vector with $|\eta| \leq \ell+1$. We can write $\eta = \eta_{0} + \tilde{\omega}$ where $\eta_{0}$ and $\tilde{\omega}$ are parallel to $\eta$ and $|\eta_{0}| \leq \ell$, $|\tilde{\omega}| = 1$. Applying the induction hypothesis to $\eta_{0}$ gives unit vectors $\omega_{1}, \ldots, \omega_{\ell}$ that add up to $\eta_{0}$. The induction step is completed by setting $\omega_{\ell+1} = \tilde{\omega}$.

To prove the theorem we choose special solutions of $\Delta v_{j} + k^{2} v_{j} = 0$ in $\Omega$ ($j=1,\dots,m+1$) having the form
\begin{equation*}
v_{j} = \rme^{\rmi \zeta_{j} \cdot x}
\end{equation*}
where $\zeta_{j} \in \mathbb{C}^{n}$ satisfy $\zeta_{j} \cdot \zeta_{j} = k^{2}$. Since $|\frac{\xi}{k}| \leq m+1$, the claim above shows that we can find unit vectors $\omega_{1}, \ldots, \omega_{m+1}$ such that
\begin{equation*}
\sum_{j=1}^{m+1} \omega_{j} = \frac{\xi}{k}.
\end{equation*}
Thus, choosing $\zeta_{j} = k \omega_{j}$, we have
\begin{equation*}
\sum_{j=1}^{m+1} \zeta_{j} = \xi.
\end{equation*}
It follows that $v_{1} \cdots v_{m+1} = \rme^{\rmi \xi \cdot x}$. Now \eref{mth_derivative_identity_estimate} shows that
\begin{equation*}
| \mathcal{F}[c](\xi) |
\leq \frac{\epsilon}{m!} \prod_{j=1}^{m+1} \| v_{j} \|_{C^{3}(\overline{\Omega})}.
\end{equation*}
The proof is completed upon observing that $\| v_{j} \|_{C^{3}(\overline{\Omega})}^2 \leq 3(1+k^{6})$ when $k \geq 1$.
\end{proof}

The assumption $|\xi| \leq (m+1) k$ ensured that we could choose solutions $v_{j} = \rme^{\rmi \zeta_{j} \cdot x}$ with $\zeta_{j}$ purely real in the proof. When $|\xi| > (m+1) k$ this will no longer be possible, and there will be a logarithmic component in the increasing stability estimate. We will next prove such an estimate for the linearized DtN map $D^{m}_{0} \Lambda_{c}$ by making a more careful choice of the vectors $\zeta_{j}$. Without loss of generality we assume that $0 \in \Omega$ and denote $D := 2 \sup_{x \in \Omega} \left| x \right|$.

\begin{thm}\label{thm_holder}
Let $D \leq 1$, $\|c\|_{C^{1}(\overline{\Omega})} \leq M_{1}$, and $k > 1$, $\epsilon < 1$, then the following estimate holds true
\begin{equation*}
\|c\|_{L^{2}(\Omega)}^{2}
\leq C k^{n+6(m+1)} \epsilon^{2} +C E^{n+6(m+1)} \epsilon
+ \frac{M_{1}^{2}}{1 + m^{2} k^{2} + E^{2}}
\end{equation*}
for the linearized system \eref{mth_equation_difference} with $E = -\ln\epsilon$ and the constant $C$ depending on the domain $\Omega$, the nonlinearity index $m$ and the dimensionality $n$.
\end{thm}

\begin{proof}

To prove the stability estimate, we shall choose the complex exponential solutions $v_{j}$ in \eref{mth_derivative_identity_estimate} carefully. Let $\xi \in \mathbb{R}^{n}$ with $\xi \neq 0$ and choose an orthonormal base $\left\{ e_{1} := \frac{\xi}{|\xi|}, e_{2}, \ldots, e_{n} \right\}$ of $\mathbb{R}^{n}$, $n \geq 2$. Let $v_{j} = \rme^{\rmi \zeta_{j} \cdot x}$ be a solution of the Helmholtz equations where the complex vectors $\zeta_{j} \in \mathbb{C}^{n}$, $j = 1,2,\ldots,m+1$ satisfy $\zeta_{j} \cdot \zeta_{j} = k^{2}$ and $\sum_{j=1}^{m+1} \zeta_{j} = \xi$.

We carry out the proof by choosing the nonlinearity index $m$ differently.

\begin{description}
  \item[Case 1: even $m$. ]
  The complex exponential solutions $v_{j} = \rme^{\rmi \zeta_{j} \cdot x}$ are constructed below by
  \begin{equation}\label{eq_thm1_CExsol1}
  \left\{~
  \eqalign{
  \zeta_{1} &= \frac{1}{m} (-k+|\xi|)e_{1} +\frac{1}{m} \sqrt{(m^{2}-1) k^{2} + 2k|\xi| - |\xi|^{2}} e_{2}, \\
  \zeta_{2} &= \frac{1}{m} (-k+|\xi|)e_{1} -\frac{1}{m} \sqrt{(m^{2}-1) k^{2} + 2k|\xi| - |\xi|^{2}} e_{2}, \\
  \zeta_{3} &= \zeta_{1}, \\
  \zeta_{4} &= \zeta_{2}, \\
  \ldots \\
  \zeta_{m-1} &=\zeta_{1}, \\
  \zeta_{m} &=\zeta_{2}, \\
  \zeta_{m+1} &= k e_{1}.
  }
  \right.
  \end{equation}
  Denote $\Xi := \sqrt{|\xi|^{2} - 2k|\xi|- (m^{2}-1) k^{2}}$. If $k \geq \frac{|\xi|}{m+1}$ and $k > 1$, then we obtain
  \begin{equation*}
  \| v_{j} \|^{2}_{C^{3}(\overline{\Omega})} \leq 3(1+k^{6}), \quad j = 1,\ldots,m+1.
  \end{equation*}
  If $k < \frac{|\xi|}{m+1}$ and $k > 1$, for $j = 1,\ldots,m$, we derive
  \begin{equation*}
  \| v_{j} \|^{2}_{C^{3}(\overline{\Omega})} \leq 3(1+k^{6}) \sup| \rme^{\rmi \zeta_{j} \cdot x}|^{2}
  \leq 3(1+k^{6}) \rme^{D\frac{\Xi}{m}}
  \end{equation*}
  and for $j = m+1$
  \begin{equation*}
  \| v_{m+1} \|^{2}_{C^{3}(\overline{\Omega})} \leq 3(1+k^{6}).
  \end{equation*}
  Recalling the identity \eref{eq_nonlinearsmallboundary_equality} and the inequality \eref{mth_derivative_identity_estimate} we have
  \begin{equation*}
  | \mathcal{F}[c](\xi) |^{2}
  = \left| \int_{\Omega} c(x) v_{1} \cdots v_{m+1} ~\rmd x \right|^{2}
  \leq \frac{\epsilon^{2}}{(m!)^{2}} \prod_{j=1}^{m+1} \| v_{j} \|^{2}_{C^{3}(\overline{\Omega})}.
  \end{equation*}
  Thus it is straightforward to obtain, for $k \geq \frac{|\xi|}{m+1}$ and $k>1$, that
  \begin{equation*}
  |\mathcal{F}[c](\xi)|^{2} \leq \frac{3^{m+1}}{(m!)^{2}} \epsilon^{2} (1+k^{6})^{m+1}
  \end{equation*}
  and for $k < \frac{|\xi|}{m+1}$, that
  \begin{equation*}
  |\mathcal{F}[c](\xi)|^{2} \leq \frac{3^{m+1}}{(m!)^{2}} \epsilon^{2} (1+k^{6})^{m+1} \rme^{D \Xi}.
  \end{equation*}
  Now we let $E := -\ln\epsilon > 0$ by assuming $\epsilon < 1$ and consider two situations such that
  \begin{description}
    \item[a)] $k > E$ (i.e.\ $\epsilon = \rme^{-E} > \rme^{-k}$) and
    \item[b)] $k \leq E$ (i.e.\ $\epsilon = \rme^{-E} \leq \rme^{-k}$).
  \end{description}
  In the situation of a), we directly obtain, with a generic constant $C := C(\Omega,m,n)$,
  \begin{equation*}
  \eqalign{
  \|c\|_{L^{2}(\Omega)}^{2}
  & = \int |\mathcal{F}[c](\xi)|^{2} ~\rmd \xi
  = \int_{k \geq \frac{|\xi|}{m+1}} |\mathcal{F}[c](\xi)|^{2} ~\rmd \xi
  + \int_{k < \frac{|\xi|}{m+1}} |\mathcal{F}[c](\xi)|^{2} ~\rmd \xi \\
  & \leq C (1+k^{6})^{m+1} (m+1)^{n} k^{n} \epsilon^{2} + \frac{M_{1}^{2}}{1+(m+1)^{2}k^{2}} \\
  & \leq C k^{n+6(m+1)} \epsilon^{2} + \frac{M_{1}^{2}}{1 + m^{2} k^{2} + E^{2}}.
  }
  \end{equation*}
  In the situation of b), we let $\rho := k + \sqrt{m^{2} k^{2}+\left(\frac{E}{D}\right)^{2}}$ such that $\sqrt{\rho^{2} - 2k\rho- (m^{2}-1)k^{2}} = \frac{E}{D}$ and split
  \begin{equation}\label{eq_thm1_c}
  \eqalign{
  \|c\|_{L^{2}(\Omega)}^{2}
  & = \int_{k \geq \frac{|\xi|}{m+1}} |\mathcal{F}[c](\xi)|^{2} ~\rmd \xi
  + \int_{k < \frac{|\xi|}{m+1} < \frac{\rho}{m+1}} |\mathcal{F}[c](\xi)|^{2} ~\rmd \xi \\
  &\quad + \int_{\rho \leq |\xi|} |\mathcal{F}[c](\xi)|^{2} ~\rmd \xi.
  }
  \end{equation}
  Meanwhile, we bound, noticing $\rho \leq (m+1)k + \frac{E}{D}$ and $k \leq E$,
  \begin{equation}\label{eq_thm1_interval}
  \eqalign{
  \int_{k < \frac{|\xi|}{m+1} < \frac{\rho}{m+1}} ~\rmd \xi
  & = \sigma_{n} \left( \rho^{n} - (m+1)^{n} k^{n} \right) \\
  & \leq \sigma_{n} \frac{E^{n}}{D^{n}} \left[ \left(1+(m+1)k\frac{D}{E}\right)^{n} - \left((m+1)k\frac{D}{E}\right)^{n} \right] \\
  & \leq \sigma_{n} \frac{E^{n}}{D^{n}} \left[ \left(1+(m+1)D\right)^{n} - \left((m+1)D\right)^{n} \right]
  }
  \end{equation}
  where $\sigma_{n}$ is the volume of an unit ball in $\mathbb{R}^{n}$.
  Then the first two terms in \eref{eq_thm1_c} can be bounded by
  \begin{equation*}
  \int_{k \geq \frac{|\xi|}{m+1}} |\mathcal{F}[c](\xi)|^{2} ~\rmd \xi
  \leq C k^{n+6(m+1)} \epsilon^{2} \leq C E^{n+6(m+1)} \epsilon^{2},
  \end{equation*}
  \begin{equation*}
  \eqalign{
  \int_{k < \frac{|\xi|}{m+1} < \frac{\rho}{m+1}} |\mathcal{F}[c](\xi)|^{2} ~\rmd \xi & \leq C k^{6(m+1)} \epsilon^{2} \rme^E \int_{k < \frac{|\xi|}{m+1} < \frac{\rho}{m+1}} ~\rmd \xi \\
  & \leq C E^{n+6(m+1)} \epsilon, 
  }
  \end{equation*}
  noticing $\int_{k < \frac{|\xi|}{m+1} < \frac{\rho}{m+1}} ~\rmd \xi \leq C E^{n}$ by \eref{eq_thm1_interval} and $k \leq E$.
  We thus obtain
  \begin{equation*}
  \eqalign{
  \|c\|_{L^{2}(\Omega)}^{2}
  & \leq C E^{n+6(m+1)} \epsilon^{2}+ C E^{n+6(m+1)} \epsilon + \frac{M_{1}^{2}}{1+m^{2}k^{2}+\frac{E^{2}}{D^{2}}} \\
  & \leq C E^{n+6(m+1)} \epsilon + \frac{M_{1}^{2}}{1+m^{2}k^{2}+E^{2}}
  }
  \end{equation*}
  since $\rho \geq \sqrt{m^{2} k^{2}+\left(\frac{E}{D}\right)^{2}}$, $\epsilon<1$ and $D \leq 1$.

  \item[Case 2: odd $m$. ]
  In this case, we could construct
  \begin{equation}\label{eq_thm1_CExsol2}
  \left\{~
  \eqalign{
  \zeta_{1} &= \frac{1}{m+1} |\xi| e_{1} + \frac{1}{m+1}\sqrt{(m+1)^{2} k^{2} -|\xi|^{2}} e_{2}, \\
  \zeta_{2} &= \frac{1}{m+1} |\xi| e_{1} - \frac{1}{m+1}\sqrt{(m+1)^{2} k^{2} -|\xi|^{2}} e_{2}, \\
  \ldots \\
  \zeta_{m} &= \zeta_{1}, \\
  \zeta_{m+1} &= \zeta_{2}.
  }
  \right.
  \end{equation}
  The analysis is similar to {\bf Case 1} by replacing $\Xi := \sqrt{|\xi|^{2}-(m+1)^{2} k^{2}}$ and $\rho := \sqrt{(m+1)^{2} k^{2} + \left(\frac{E}{D}\right)^{2}}$. If $k > E$, we obtain
  \begin{equation*}
  \eqalign{
  \|c\|_{L^{2}(\Omega)}^{2}
  &= \int_{k \geq \frac{|\xi|}{m+1}} |\mathcal{F}[c](\xi)|^{2} ~\rmd \xi
  + \int_{k < \frac{|\xi|}{m+1}} |\mathcal{F}[c](\xi)|^{2} ~\rmd \xi \\
  & \leq C (1+k^{6})^{m+1} (m+1)^{n} k^{n} \epsilon^{2} + \frac{M_{1}^{2}}{1+(m+1)^{2}k^{2}} \\
  & \leq C k^{n+6(m+1)} \epsilon^{2} + \frac{M_{1}^{2}}{1 + m^{2} k^{2} + E^{2}}.
  }
  \end{equation*}
  If $k \leq E$ and $D \leq 1$, we have
  \begin{equation*}
  \eqalign{
  \|c\|_{L^{2}(\Omega)}^{2}
  &= \int_{k \geq \frac{|\xi|}{m+1}} |\mathcal{F}[c](\xi)|^{2} ~\rmd \xi
  + \int_{k < \frac{|\xi|}{m+1} < \frac{\rho}{m+1}} |\mathcal{F}[c](\xi)|^{2} ~\rmd \xi \\
  &\quad + \int_{\rho \leq |\xi|} |\mathcal{F}[c](\xi)|^{2} ~\rmd \xi \\
  & \leq C E^{n+6(m+1)} \epsilon^{2} + C E^{n+6(m+1)} \epsilon + \frac{M_{1}^{2}}{1+(m+1)^{2}k^{2}+E^{2}}.
  }
  \end{equation*}
\end{description}
\end{proof}

\begin{rem}
We shall mention that treatment of the identity \eref{eq_nonlinearsmallboundary_equality} in current work is quite different from that in \cite{LLLS_I}. More precisely, in \cite{LLLS_I}, the authors consider an inverse problem for elliptic equations where $v_{j}$ are solutions of Laplace equations. Since any constant is a trivial solution there, the uniqueness in \cite{LLLS_I} can be obtained based on the classic arguments in \cite{C1980}. On the other hand, in current work, $v_{j}$ represent the solutions of Helmholtz equations and we have to choose them very carefully as shown in the above proof.
\end{rem}

Despite the profound theoretical justification by the linearized DtN map $D^{m}_{0} \Lambda_{c}$, it is somehow not easy to approximate such a linearized DtN map numerically which will be shown in \Sref{se4.3}. In particular, the small boundary data yields a solution with small values which is easily contaminated by numerical differentiation error. To further study the linearized inverse Schr\"{o}dinger potential problem of \eref{eq_HelmholtzEqmain}, it is worthwhile to consider the linearized DtN map $\Lambda'_{c}$ corresponding to the case where $c$ is small compared to the wavenumber $k$. In particular, we prove the uniqueness for the linearized inverse Schr\"{o}dinger potential problem with an arbitrary power type nonlinearity term below given the linearized DtN map $\Lambda'_{c}$ at a fixed wavenumber $k > 0$.

\begin{thm}\label{thm_unique_m}
Let $c_{1}$ and $c_{2}$ be two functions in $L^{\infty}(\Omega)$. If the two linearized DtN maps in \eref{eq_LinearDtN} obey $\Lambda'_{c_{1}} = \Lambda'_{c_{2}}$, then $c_{1} = c_{2}$ in $\Omega$.
\end{thm}

\begin{proof}
The proof is similar to the seminal work by Calder\'{o}n \cite{C1980} but one needs to choose appropriate complex exponential solutions. To this end, we let $\xi \in \mathbb{R}^{n}$ with $\xi \neq 0$ and choose an orthonormal base $\left\{ e_{1} := \frac{\xi}{|\xi|}, e_{2}, \ldots, e_{n} \right\}$ of $\mathbb{R}^{n}$, $n \geq 2$. Then we can define the following vectors $\mu_{\ell} \in \mathbb{C}^{n}$, $\ell = 1,2$ such that
\begin{equation}\label{eq_exposol_first}
\left\{~
\eqalign{
\mu_{1} & = \frac{+(m^{2}-1) k^{2} + |\xi|^{2}}{2m|\xi|} e_{1} - \frac{\sqrt{-(m^{2}-1)^{2} k^{4} + 2(m^{2}+1) k^{2} |\xi|^{2} -|\xi|^{4}}}{2m|\xi|} e_{2}, \\
\mu_{2} & = \frac{-(m^{2}-1) k^{2} + |\xi|^{2}}{2 |\xi|} e_{1} + \frac{\sqrt{-(m^{2}-1)^{2} k^{4} + 2(m^{2}+1) k^{2} |\xi|^{2} -|\xi|^{4}}}{2 |\xi|} e_{2}.
}
\right.
\end{equation}
We end the proof by assigning, in \eref{eq_calderonILX},
\begin{equation}\label{eq_exposol_m}
u_{0}(x) = \rme^{\rmi \mu_{1} \cdot x}, \qquad
\varphi(x) = \rme^{\rmi \mu_{2} \cdot x},
\end{equation}
such that $u_{0}^{m}(x) \varphi(x) = \rme^{\rmi \xi \cdot x}$.
\end{proof}

\begin{rem}\label{rem_uniqueComSol}
For $m = 1$, namely, when the power-type nonlinearity term reduces to a linear one, the complex exponential solutions in \eref{eq_exposol_m} are exactly those solutions used in \cite{ILX2020}. Nevertheless, it is somehow disappointing that when $m \geq 2$, \eref{eq_exposol_m} does not easily give a stability estimate.
In particular, the failure is exactly induced by the behavior of the complex vectors in \eref{eq_exposol_first}. It is easy to verify that
\begin{equation*}
\eqalign{
&-(m^{2}-1)^{2} k^{4} + 2(m^{2}+1) k^{2} |\xi|^{2} - |\xi|^{4} \\
&\quad = - \big(|\xi|+(m+1)k\big) \big(|\xi|+(m-1)k\big) \big(|\xi|-(m-1)k\big) \big(|\xi|-(m+1)k\big)
\geq 0,
}
\end{equation*}
for $|\xi| \in [(m-1)k,(m+1)k]$. If $|\xi| \in \big(0,(m-1)k\big)$ or $|\xi| > (m+1)k$, the complex exponential solution $u_{0}(x) = \rme^{\rmi \mu_{1} \cdot x}$ or $\varphi(x) = \rme^{\rmi \mu_{2} \cdot x}$ blows up at $e_{2}$ (or $-e_{2}$) direction when $|x|$ increases.
\end{rem}

Though it is not straightforward to obtain an increasing stability by the linearized DtN map $\Lambda'_{c}$ for an arbitrary choice of the nonlinearity index $m$, as shown in Remark \ref{rem_uniqueComSol}, we could still stably reconstruct the Fourier coefficients of the unknown potential function $c(x)$ within an interval given any fixed wavenumber $k > 0$. This observation allows us to design a reconstruction algorithm if the linearized DtN map $\Lambda'_{c}$ of multiple wavenumbers are provided. We will discuss this in \Sref{se4.2}.

On the other hand, if one chooses a specific power-type nonlinearity term, for instance a quadratic one with $m = 2$, we could regain the increasing stability by calibrating the identity \eref{eq_calderonILX} carefully. This result will be given in the coming \Sref{se3}.


\section{Linearized inverse Schr\"{o}dinger potential problem with a quadratic nonlinearity term}\label{se3}

To obtain a stability estimate of the linearized inverse Schr\"{o}dinger potential problem with a power type nonlinearity term by the linearized DtN map $\Lambda'_{c}$, the construction of the complex exponential solutions is essential and the standard approach in \Sref{se2} fails in view of the discussion in Remark \ref{rem_uniqueComSol}. To successfully prove the stability estimate, we may have to treat the nonlinearity term separately and the linearized inverse Schr\"{o}dinger potential problem with a quadratic nonlinearity term ($m = 2$) will be extensively investigated in current section. For the situation of a general nonlinearity index $m > 2$, we consider it as a future work and will report the result elsewhere.

To proceed further, we are inspired by the idea of small boundary data discussed above and consider three solutions of \eref{eq_HelmholtzEqmain} which are denoted by $u$, $v$ and $w$ with appropriate Dirichlet boundary conditions. By assuming that the potential function $c(x)$ is small or the squared wavenumber $k^{2}$ is sufficiently large, and recalling the asymptotical expansion of these solutions as in \eref{eq_asymILX}, we have
\begin{equation*}
\eqalign{
u &= u_{0} + u_{1} + \ldots, \\
v &= v_{0} + v_{1} + \ldots, \\
w &= w_{0} + w_{1} + \ldots,
}
\end{equation*}
where the remaining ``$\ldots$" are the higher order terms of these solutions. In fact, we can obtain that
\begin{equation}\label{eq_quadratic_uv}
\left\{~
\eqalign{
\Delta u_{0} + k^{2} u_{0} = 0, \quad
\Delta u_{1} + k^{2} u_{1} = c(x) u_{0}^{2} &\quad \textrm{in\ } \Omega, \\
\Delta v_{0} + k^{2} v_{0} = 0, ~\quad
\Delta v_{1} + k^{2} v_{1} = c(x) v_{0}^{2} &\quad \textrm{in\ } \Omega,
}
\right.
\end{equation}
for $u$, $v$ and for $w$,
\begin{equation}\label{eq_quadratic_w}
\Delta w_{0} + k^{2} w_{0} = 0, \quad
\Delta w_{1} + k^{2} w_{1} = c(x) w_{0}^{2} \quad \textrm{in\ } \Omega.
\end{equation}
The linearized DtN map $\Lambda'_{c}$ can be defined accordingly for these three solutions as in \eref{eq_LinearDtN}.

Denoting the Dirichlet boundary conditions of $u_{0}$, $v_{0}$ by $u_{0}|_{\partial\Omega}$ and $v_{0}|_{\partial\Omega}$, we define the boundary condition of $w_{0}$ by
\begin{equation*}
w_{0}|_{\partial\Omega} := u_{0}|_{\partial\Omega} + v_{0}|_{\partial\Omega}.
\end{equation*}
By the linearity of the Helmholtz equation, we know
\begin{equation*}
w_{0} = u_{0} + v_{0} \quad \textrm{in\ } \Omega.
\end{equation*}

We now take a close look of the asymptotical expansion of $w = w_{0} + w_{1} + \ldots$ in \eref{eq_quadratic_w} and choose $\varphi$ to be another solution of the Helmholtz equation $\Delta \varphi + k^{2} \varphi = 0$ in $\Omega$, then, while $w_{1}|_{\partial\Omega} = 0$, we have
\begin{equation*}
\int_{\Omega} c(x) w_{0}^{2} \varphi ~\rmd x
= \int_{\partial\Omega} \partial_{\nu} w_{1} \varphi ~\rmd S.
\end{equation*}
Noticing that $w_{0}^{2} = u_{0}^{2} + v_{0}^{2} + 2 u_{0} v_{0}$ in $\Omega$, we thus obtain
\begin{equation*}
2 \int_{\Omega} c(x) u_{0} v_{0} \varphi ~\rmd x
= \int_{\Omega} c(x) w_{0}^{2} \varphi ~\rmd x
- \left( \int_{\Omega} c(x) u_{0}^{2} \varphi ~\rmd x
+ \int_{\Omega} c(x) v_{0}^{2} \varphi ~\rmd x \right).
\end{equation*}
Recalling the asymptotical expansion of $u$ and $v$, as $u_{1}|_{\partial\Omega} = 0$ and $v_{1}|_{\partial\Omega} = 0$, we derive
\begin{equation}\label{eq_Calderonformula}
2 \int_{\Omega} c(x) u_{0} v_{0} \varphi ~\rmd x
= \int_{\partial\Omega} \partial_{\nu} w_{1} \varphi ~\rmd S
- \left( \int_{\partial\Omega} \partial_{\nu} u_{1} \varphi ~\rmd S
+ \int_{\partial\Omega} \partial_{\nu} v_{1} \varphi ~\rmd S \right).
\end{equation}
The identity \eref{eq_Calderonformula} then allows us to carry out the stability estimate and reconstruction algorithm of the linearized inverse Schr\"{o}dinger potential problem with a quadratic nonlinearity term.

Similar to Theorem \ref{thm_holder}, we again assume $0 \in \Omega$, $D = 2 \sup_{x\in\Omega} |x|$ and denote the same variable $\epsilon : = \sup_{\|\tilde{g}_{0}\|_{C^{2}(\partial\Omega)} = 1} \|\Lambda'_{c} \tilde{g}_{0}\|_{L^{\frac{3}{2}}(\partial\Omega)}$ to be the operator norm of $\Lambda'_{c}$ defined in \eref{eq_LinearDtN}. Then for any $g_{0} \in C^{2}(\partial\Omega)$, $G: =\|g_{0}\|_{C^{2}(\partial\Omega)}$ and define $\tilde{u}_{0}$ be the solution of
\begin{equation*}
\left\{~
\eqalign{
\Delta \tilde{u}_{0} + k^{2} \tilde{u}_{0} = 0 &\quad \textrm{in\ } \Omega, \\
\tilde{u}_{0} = g_{0}/G &\quad \textrm{on\ } \partial\Omega,
}
\right.
\end{equation*}
we have the solution $u_{1}$ of \eref{eq_quadratic_u1} by
\begin{equation*}
\Delta u_{1} + k^{2} u_{1} = c(x) G^2 \tilde{u}_{0}^{2}  \quad \textrm{in\ } \Omega,
\end{equation*}
and $\|\partial_{\nu} u_{1}\|_{L^{\frac{3}{2}}(\partial\Omega)} = \|\Lambda'_{c} g_{0}\|_{L^{\frac{3}{2}}(\partial\Omega)} \leq \epsilon \|g_{0}\|_{C^{2}(\partial\Omega)}^{2}$ consequently. The main stability estimate is presented below.

\begin{thm}\label{thm_quadratic}
Let $D \leq 1$, $\|c\|_{H^{1}(\Omega)} \leq M_{1}$, and $k > 1$, $\epsilon < 1$, then the following estimate holds true
\begin{equation*}
\|c\|_{L^{2}(\Omega)}^{2}
\leq C \left( k^{n+8} + E^{n+8} \right) \epsilon^{2}
+ C E^{n+8} \epsilon
+ \frac{M_{1}^{2}}{1 + 4 k^{2} + E^{2}}
\end{equation*}
for the linearized system \eref{eq_quadratic_uv}, \eref{eq_quadratic_w} with $E = -\ln\epsilon$ and the constant $C$ depending on the domain $\Omega$ and the dimensionality $n$.
\end{thm}

\begin{proof}

Let $\xi \in \mathbb{R}^{n}$ with $\xi \neq 0$ and choose an orthonormal base $\left\{ e_{1} := \frac{\xi}{|\xi|}, e_{2}, \ldots, e_{n} \right\}$ of $\mathbb{R}^{n}$, $n \geq 2$. Then we can choose the following $\zeta_{\ell} \in \mathbb{C}^{n}$, $\ell = 1,2,3$ such that
\begin{equation*}
\left\{~
\eqalign{
\zeta_{1} & = \frac{1}{2} (-k+|\xi|) e_{1} - \frac{1}{2} \sqrt{3k^{2}+2k|\xi|-|\xi|^{2}} e_{2}, \\
\zeta_{2} & = \frac{1}{2} (-k+|\xi|) e_{1} + \frac{1}{2} \sqrt{3k^{2}+2k|\xi|-|\xi|^{2}} e_{2}, \\
\zeta_{3} & = k e_{1}.
}
\right.
\end{equation*}
We assign
\begin{equation}\label{eq_staproof_ComSol}
u_{0}(x) = \rme^{\rmi \zeta_{1} \cdot x}, \qquad
v_{0}(x) = \rme^{\rmi \zeta_{2} \cdot x}, \qquad
\varphi(x) = \rme^{\rmi \zeta_{3} \cdot x}.
\end{equation}
Then
\begin{equation*}
u_{0} v_{0} \varphi = \rme^{\rmi \xi \cdot x}
\end{equation*}
and the identity \eref{eq_Calderonformula} yields
\begin{equation*}
2 \mathcal{F}[c](\xi) = 2 \int_{\Omega} c(x) \rme^{\rmi \xi \cdot x} ~\rmd x
= \int_{\partial\Omega} \partial_{\nu} w_{1} \varphi ~\rmd S
- \left( \int_{\partial\Omega} \partial_{\nu} u_{1} \varphi ~\rmd S
+ \int_{\partial\Omega} \partial_{\nu} v_{1} \varphi ~\rmd S \right).
\end{equation*}
Hence we obtain, since $w_{0} = u_{0} + v_{0}$ in $\Omega$,
\begin{equation*}
\eqalign{
|\mathcal{F}[c](\xi)|^{2}
& \leq \frac{1}{4} \left(
\|\partial_{\nu} w_{1}\|_{L^{\frac{3}{2}}(\partial\Omega)}^{2}
+ \|\partial_{\nu} u_{1}\|_{L^{\frac{3}{2}}(\partial\Omega)}^{2}
+ \|\partial_{\nu} v_{1}\|_{L^{\frac{3}{2}}(\partial\Omega)}^{2}
\right) \|\varphi\|_{L^{3}(\partial\Omega)}^{2} \\
& \leq \frac{1}{4} \epsilon^{2} \left(
\|w_{0}|_{\partial\Omega}\|_{C^{2}(\partial\Omega)}^{4}
+ \|u_{0}|_{\partial\Omega}\|_{C^{2}(\partial\Omega)}^{4}
+ \|v_{0}|_{\partial\Omega}\|_{C^{2}(\partial\Omega)}^{4}
\right) \|\varphi\|_{L^{3}(\partial\Omega)}^{2} \\
& \leq C \epsilon^{2} \left(
\|w_{0}\|_{C^{2}(\overline{\Omega})}^{4}
+ \|u_{0}\|_{C^{2}(\overline{\Omega})}^{4}
+ \|v_{0}\|_{C^{2}(\overline{\Omega})}^{4}
\right) \|\varphi\|_{L^{\infty}(\Omega)}^{2} \\
& \leq C \epsilon^{2} \left(
\|u_{0}\|_{C^{2}(\overline{\Omega})}^{4}
+ \|v_{0}\|_{C^{2}(\overline{\Omega})}^{4}
\right)
}
\end{equation*}
with a generic constant $C$ depending on the domain $\Omega$, and $\|\varphi\|_{L^{\infty}(\Omega)}^{2} = \|\rme^{\rmi k e_{1} \cdot x}\|_{L^{\infty}(\Omega)}^{2} \leq 1$.

Noticing the fact that $|\zeta_{\ell}|^{2} = k^{2}$, $\ell = 1,2,3$, we thus obtain, if $k \geq \frac{|\xi|}{3}$,
\begin{equation*}
\|u_{0}\|_{C^{2}(\overline{\Omega})}^{4}
= \|v_{0}\|_{C^{2}(\overline{\Omega})}^{4}
\leq C \left(1+k^{8}\right).
\end{equation*}
Then, there holds
\begin{equation*}
|\mathcal{F}[c](\xi)|^{2}
\leq C \epsilon^{2} \left(1+k^{8}\right), \qquad \textrm{for\ } k \geq \frac{|\xi|}{3}.
\end{equation*}

If $k < \frac{|\xi|}{3}$, by denoting $\Xi := \sqrt{|\xi|^{2}-2k|\xi|-3k^{2}}$ we then derive the following bounds
\begin{equation*}
\|u_{0}\|_{C^{2}(\overline{\Omega})}^{4}
= \|v_{0}\|_{C^{2}(\overline{\Omega})}^{4}
\leq C \left(1+k^{8}\right) \|\rme^{\rmi \zeta_{1} \cdot x}\|_{L^{\infty}(\Omega)}^{4}
\leq C \left(1+k^{8}\right) \rme^{D\Xi}.
\end{equation*}
Consequently, we derive
\begin{equation*}
|\mathcal{F}[c](\xi)|^{2}
\leq C \epsilon^{2} \left(1+k^{8}\right) \rme^{D\Xi}, \qquad \textrm{for\ } k < \frac{|\xi|}{3}.
\end{equation*}

Let $E := -\ln\epsilon > 0$ and $k > 1$, $\epsilon < 1$, we again consider two cases
\begin{itemize}
\item[a)] $k > E$ (i.e.\ $\epsilon = \rme^{-E} > \rme^{-k}$), and
\item[b)] $k \leq E$ (i.e.\ $\epsilon = \rme^{-E} \leq \rme^{-k}$).
\end{itemize}

In the case a), we have
\begin{equation*}
\eqalign{
\|c\|_{L^{2}(\Omega)}^{2}
& = \int |\mathcal{F}[c](\xi)|^{2} ~\rmd \xi
= \int_{k \geq \frac{|\xi|}{3}} |\mathcal{F}[c](\xi)|^{2} ~\rmd \xi
+ \int_{k < \frac{|\xi|}{3}} |\mathcal{F}[c](\xi)|^{2} ~\rmd \xi \\
& \leq C \epsilon^{2} \left(1+k^{8}\right) \sigma_{n} (3k)^{n} + \frac{M_{1}^{2}}{1+(3k)^{2}} \\
& \leq C k^{n+8} \epsilon^{2} + \frac{M_{1}^{2}}{1+8k^{2}+E^{2}}
}
\end{equation*}
where $\sigma_{n}$ is the volume of an unit ball in $\mathbb{R}^{n}$, and the constant $C$ depends on the domain $\Omega$ and the dimensionality $n$.

In the case b), we let $\rho := k + \sqrt{4k^{2}+\left(\frac{E}{D}\right)^{2}}$ such that $\sqrt{\rho^{2} - 2k \rho - 3k^{2}} = \frac{E}{D}$ and split
\begin{equation}\label{eq_proofthm31_c}
\eqalign{
\|c\|_{L^{2}(\Omega)}^{2}
& = \int_{k \geq \frac{|\xi|}{3}} |\mathcal{F}[c](\xi)|^{2} ~\rmd \xi
+ \int_{k < \frac{|\xi|}{3} < \frac{\rho}{3}} |\mathcal{F}[c](\xi)|^{2} ~\rmd \xi \\
&\quad + \int_{\rho \leq |\xi|} |\mathcal{F}[c](\xi)|^{2} ~\rmd \xi.
}
\end{equation}
The first term in the right-hand side of \eref{eq_proofthm31_c} can be bounded by
\begin{equation*}
\int_{k \geq \frac{|\xi|}{3}} |\mathcal{F}[c](\xi)|^{2} ~\rmd \xi
\leq C k^{n+8} \epsilon^{2}
\leq C E^{n+8} \epsilon^{2},
\end{equation*}
noticing $k \leq E$.

We focus on the second term in the right-hand side of \eref{eq_proofthm31_c} and estimate
\begin{equation*}
\int_{k < \frac{|\xi|}{3} < \frac{\rho}{3}} |\mathcal{F}[c](\xi)|^{2} ~\rmd \xi
\leq C \epsilon^{2} k^{8} \left( \int_{k < \frac{|\xi|}{3} < \frac{\rho}{3}} \rme^{D\Xi} ~\rmd \xi \right)
\leq C \epsilon k^{8} \left( \int_{k < \frac{|\xi|}{3} < \frac{\rho}{3}} ~\rmd \xi \right)
\end{equation*}
since $\rme^{D\Xi} \leq e^{E} = \epsilon^{-1}$ when $k < \frac{|\xi|}{3} < \frac{\rho}{3}$.
Meanwhile, noticing $\rho \leq 3k + \frac{E}{D}$ and $k \leq E$, we bound
\begin{equation*}
\eqalign{
\int_{k < \frac{|\xi|}{3} < \frac{\rho}{3}} ~\rmd \xi
& = \sigma_{n} \left( \rho^{n} - (3k)^{n} \right) \\
& \leq \sigma_{n} \frac{E^{n}}{D^{n}} \left[ \left(1+3k\frac{D}{E}\right)^{n} - \left(3k\frac{D}{E}\right)^{n} \right] \\
& \leq \sigma_{n} \frac{E^{n}}{D^{n}} \left[ \left(1+3D\right)^{n} - \left(3D\right)^{n} \right]
}
\end{equation*}
where $\sigma_{n}$ is the volume of an unit ball in $\mathbb{R}^{n}$. The above inequalities yields
\begin{equation*}
\int_{k < \frac{|\xi|}{3} < \frac{\rho}{3}} |\mathcal{F}[c](\xi)|^{2} ~\rmd \xi
\leq C \epsilon k^{8} \left( \sigma_{n} \frac{E^{n}}{D^{n}} \left[ \left(1+3D\right)^{n} - \left(3D\right)^{n} \right] \right)
\leq C E^{n+8} \epsilon.
\end{equation*}
Furthermore, since $\rho > \sqrt{4k^{2}+\left(\frac{E}{D}\right)^{2}}$, we bound the third term in the right-hand side of \eref{eq_proofthm31_c} by
\begin{equation*}
\int_{\rho \leq |\xi|} |\mathcal{F}[c](\xi)|^{2} ~\rmd \xi
\leq \frac{M_{1}^{2}}{1+\rho^{2}}
\leq \frac{M_{1}^{2}}{1+4k^{2}+\frac{E^{2}}{D^{2}}}.
\end{equation*}

We thus prove for both cases the proposed bound.

\end{proof}

\begin{rem}\label{rem_se3}
The stability estimate in above Theorem \ref{thm_quadratic}, if $k$ is sufficiently large, is similar as in \cite[Theorem 2.1]{ILX2020} where a linear elliptic equation is investigated ibid, i.e.\
\begin{equation*}
\left\{~
\eqalign{
\Delta u + k^{2} u - c(x) u = 0 &\quad \textrm{in\ } \Omega, \\
u = g_{0} &\quad \textrm{on\ } \partial\Omega.
}
\right.
\end{equation*}
A clear numerical evidence will be provided in \Sref{se4} and one can stably recover the Fourier coefficients with $|\xi| \leq 3k$ whereas in \cite{ILX2020} one can only recover those with $|\xi| \leq 2k$. Such gain highly depends on the sophisticatedly selected complex exponential functions and the modified identity \eref{eq_Calderonformula} considered above. It can be viewed as the advantage of the quadratic nonlinearity term when we solve the linearized inverse problems \eref{eq_HelmholtzEqmain} with $m = 2$.
\end{rem}


\section{Reconstruction algorithm and numerical examples}\label{se4}

In this section, we provide two reconstruction algorithms stably recovering the unknown potential function by the linearized DtN map $\Lambda'_{c}$ and a vanilla reconstruction algorithm by the linearized DtN map $D^{m}_{0} \Lambda_{c}$. In view of the quadratic nonlinearity term, we rely on the theoretical discussion in \Sref{se3} and deliver the first algorithm where boundary measurements of a single (large) wavenumber could offer a high resolution. Meanwhile, the second algorithm focuses on the high-order nonlinearity term discussed in \Sref{se2} and the linearized DtN map $\Lambda'_{c}$ of multiple wavenumbers is included to recover sufficiently many Fourier coefficients of the unknown potential function. Finally a vanilla reconstruction algorithm by the (approximated) linearized DtN map $D^{m}_{0} \Lambda_{c}$ is presented to verify the feasibility of the proposed linearization which, to the best of our knowledge, is the first attempt to realize the linearized DtN map $D^{m}_{0} \Lambda_{c}$ numerically.

\reviewA{We shall emphasize that by implementing the linearized DtN maps $\Lambda'_{c}$ and $D^{m}_{0} \Lambda_{c}$ the range of stably reconstructed Fourier mode $\mathcal{F}[c](\xi)$ is expanded to $|\xi| \leq (m+1)k$ as shown in Theorem \ref{thm_holder} with an arbitrary finite integer $m$ and Theorem \ref{thm_quadratic} with $m = 2$. This truncated value $(m+1)k$ can be viewed as a regularization for stably recovering the unknown potential function $c(x)$.}


\subsection{Reconstruction algorithm by $\Lambda'_{c}$ for a quadratic nonlinearity term}\label{se4.1}

Noticing that the linearized DtN map $\Lambda'_{c}$ can be numerically approximated, see \cite[Eq.(4.3)]{ILX2020}, we present the first reconstruction algorithm based on the identity \eref{eq_Calderonformula}. As an illustration, we focus on the two-dimensional space $n = 2$.

By selecting the complex exponential solutions \eref{eq_staproof_ComSol} in the proof of Theorem \ref{thm_quadratic}, we know that the left-hand side of \eref{eq_Calderonformula} reflects a Fourier coefficient of the potential function $c(x)$. Then by choosing $\xi \in \mathbb{R}^{n}$ and recalling Remark \ref{rem_se3}, we aim to recovering all the Fourier coefficients $\mathcal{F}[c](\xi)$ of the potential function $c(x)$ satisfying $|\xi| \leq 3k$. The larger wavenumber $k$, the more Fourier coefficients can be recovered.

To further address the reconstruction algorithm, we need the following discrete sets of lengths and angles of the vectors in the phase space. The discrete and finite length set is defined by
\begin{equation*}
\{\kappa_{i}\}_{i=1}^{I} \subset (0, L k \,]
\quad \textrm{for any fixed } k.
\end{equation*}
Here we choose $L \geq 3$ and $L k$ is the maximum length of the vector $\xi$. Two angle sets are defined by
\begin{equation*}
\{\hat{y}_{s}\}_{s=1}^{S} \subset \mathbb{S}^{n-1}
\quad \textrm{and} \quad \{\hat{z}_{s}\}_{s=1}^{S} \subset \mathbb{S}^{n-1},
\end{equation*}
which satisfy $\hat{y}_{s} \cdot \hat{z}_{s} = 0$.

The vector $\xi^{\langle i;s \rangle} := \kappa_{i} \hat{y}_{s} $ and following vectors $\zeta_{\ell}^{\langle i;s \rangle} \in \mathbb{C}^{n}$, $\ell = 1,2,3$ are chosen
\begin{equation*}
\left\{~
\eqalign{
\zeta_{1}^{\langle i;s \rangle} &:= \frac{1}{2} (-k+\kappa_{i}) \hat{y}_{s} - \frac{1}{2} \sqrt{3k^{2}+2k\kappa_{i}-\kappa_{i}^{2}} \hat{z}_{s}, \\
\zeta_{2}^{\langle i;s \rangle} &:= \frac{1}{2} (-k+\kappa_{i}) \hat{y}_{s} + \frac{1}{2} \sqrt{3k^{2}+2k\kappa_{i}-\kappa_{i}^{2}} \hat{z}_{s}, \\
\zeta_{3}^{\langle i;s \rangle} &:= k \hat{y}_{s},
}
\right.
\end{equation*}
which further assign to the complex exponential solution as in \eref{eq_staproof_ComSol} below
\begin{equation*}
u_{0}(x) = \rme^{\rmi \zeta_{1}^{\langle i;s \rangle} \cdot x}, \qquad
v_{0}(x) = \rme^{\rmi \zeta_{2}^{\langle i;s \rangle} \cdot x}, \qquad
\varphi(x) = \rme^{\rmi \zeta_{3}^{\langle i;s \rangle} \cdot x}
\end{equation*}
for $i = 1,2,\cdots,I$ and $s = 1,2,\cdots,S$. Here, the superscript notation $\cdot^{\langle i;s \rangle}$ will be referred to a vector $\xi^{\langle i;s \rangle}$ with the $i$th length $\kappa_{i}$ and the $s$th angle $\hat{y}_{s}$. Finally, for the inverse Fourier transform, a numerical quadrature rule can be constructed by a suitable choice of the weights $\sigma^{\langle i;s \rangle}$ according to these points $\xi^{\langle i;s \rangle}$.

We summarize our reconstruction algorithm below, which is similar to that in \cite{ILX2020} but one has to solve the nonlinear Schr\"{o}dinger potential problem three times at each iteration because of the quadratic nonlinearity term.
\vspace{10pt}
\hrule\hrule
\vspace{8pt}
{\parindent 0pt \bf Algorithm 1: Reconstruction Algorithm for the Linearized Schr\"{o}dinger Potential Problem, the quadratic nonlinearity term} %
\vspace{5pt}
\hrule
\vspace{8pt}
{\parindent 0pt \bf Input:} %
$k$, %
$\{\kappa_{i}\}_{i=1}^{I}$, %
$\{\hat{y}_{s}\}_{s=1}^{S}$, $\{\hat{z}_{s}\}_{s=1}^{S}$ and %
$\{\sigma^{\langle i;s \rangle}\}$; \\[5pt]%
{\parindent 0pt \bf Output:} %
Approximated Potential $c^{\langle I+1;1 \rangle}$. \\[-15pt]%
\begin{enumerate}
  \item[1:] $\,$ Set $c^{\langle 1;1 \rangle} := 0$; %
  \item[2:] $\,$ {\bf For} $i = 1,2,\dots,I$ (length~updating) %
  \item[3:] $\,$ \quad {\bf For} $s = 1,2,\dots,S$ (angle~updating) %
  \item[4:] $\,$ \quad \quad Choose $u_{0} := \exp \{ \rmi \zeta_{1}^{\langle i;s \rangle} \cdot x \}$, $v_{0} := \exp \{ \rmi \zeta_{2}^{\langle i;s \rangle} \cdot x \}$ and $w_{0} := u_{0} + v_{0}$; %
  \item[5:] $\,$ \quad \quad Measure the Neumann boundary data $\partial_{\nu} u$, $\partial_{\nu} v$, $\partial_{\nu} w$ of the forward problem \eref{eq_HelmholtzEqmain} %
  \item[] $\,$ \quad \qquad while the Dirichlet boundary data $u_{0}|_{\partial\Omega}$, $v_{0}|_{\partial\Omega}$, $w_{0}|_{\partial\Omega}$ are given; %
  \item[6:] $\,$ \quad \quad Calculate the approximated linearized Neumann boundary data %
  \item[] $\,$ \quad \qquad $g_{u}^{\,\prime} := (\partial_{\nu} u - \partial_{\nu} u_{0})|_{\partial\Omega}$, \quad $g_{v}^{\,\prime} := (\partial_{\nu} v - \partial_{\nu} v_{0})|_{\partial\Omega}$, \quad $g_{w}^{\,\prime} := (\partial_{\nu} w - \partial_{\nu} w_{0})|_{\partial\Omega}$; %
  \item[7:] $\,$ \quad \quad Choose $\varphi := \exp \{ \rmi \zeta_{3}^{\langle i;s \rangle} \cdot x \}$ and $\gamma := \big[ u_{0} v_{0} \varphi \big]^{-1} = \exp \{ -\rmi \xi^{\langle i;s \rangle} \cdot x \}$; %
  \item[8:] $\,$ \quad \quad Compute $\mathcal{F}[c](\xi^{\langle i;s \rangle}) \approx \frac{1}{2} \int_{\partial\Omega} ( g_{w}^{\,\prime} - g_{u}^{\,\prime} - g_{v}^{\,\prime} ) \, \varphi \,\rmd S$; %
  \item[9:] $\,$ \quad \quad Update $c^{\langle i;s+1 \rangle} := c^{\langle i;s \rangle} + \mathcal{F}[c](\xi^{\langle i;s \rangle}) \, \gamma \sigma^{\langle i;s \rangle}$, \quad if $\kappa_{i} \leq 3k$; %
  \item[10:] $\,$ \quad {\bf End}; %
  \item[11:] $\,$ \quad Set $c^{\langle i+1;1 \rangle} := c^{\langle i;S+1 \rangle}$; %
  \item[12:] $\,$ {\bf End}. %
\end{enumerate}
\vspace{8pt}
\hrule\hrule
\vspace{10pt}

In fact, the linearized Neumann boundary data $\partial_{\nu} w_{1}$ depends on the unknown potential function $c(x)$ referring to \eref{eq_quadratic_w}. As mentioned in \cite[Eq.(4.3)]{ILX2020}, we utilize
\begin{equation}\label{eq_se4linearizedNeumann}
g_{w}^{\,\prime} := (\partial_{\nu} w - \partial_{\nu} w_{0})|_{\partial\Omega}
\end{equation}
to approximate the non-measurable data $\partial_{\nu} w_{1}|_{\partial\Omega}$. The similar approximation $g_{u}^{\,\prime}$ and $g_{v}^{\,\prime}$ are employed for the linearized Neumann data $\partial_{\nu} u_{1}|_{\partial\Omega}$ and $\partial_{\nu} v_{1}|_{\partial\Omega}$, respectively.

As one can observe, the computational cost of \textbf{Algorithm 1} is quite high because of the nonlinearity term in the forward problem, see e.g. \cite{FT2005, WZ2018, XB2010, YL2017}. In particular in Steps 5-6 of \textbf{Algorithm 1}, we must solve the nonlinear elliptic equation \eref{eq_HelmholtzEqmain} three times in order to derive their Neumann traces which are necessary to compute the Fourier coefficient in Step 8.

\begin{figure}[htbp]
\centering
\includegraphics[width=0.5\textwidth]{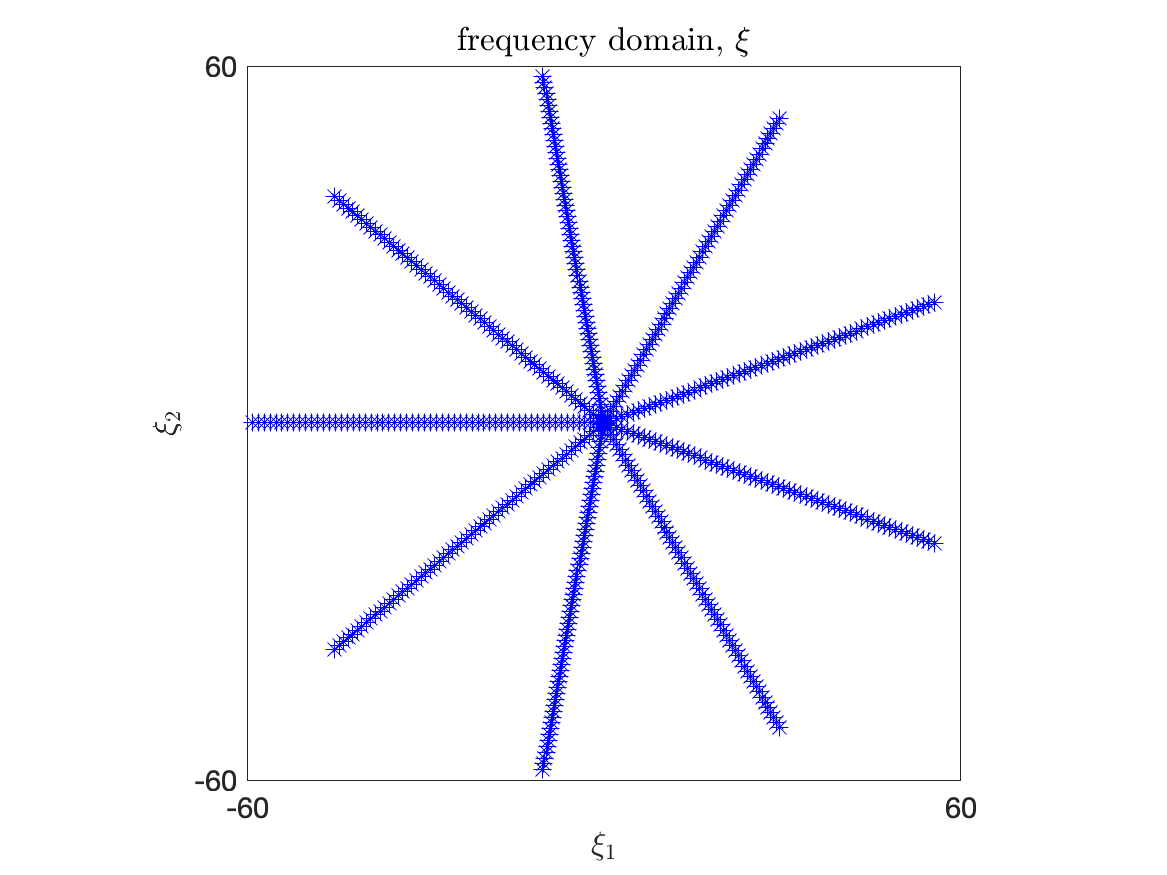}\\
\caption{The sampling points $\xi = (\xi_{1},\xi_{2})$ in frequency domain.}
\label{fig:Xi}
\end{figure}

\begin{figure}[htbp]
\centering
\,\hfill \textbf{Quadratic case}: \hfill\,\\
\,\hfill \textbf{(a)} $k = 5$ \hfill\,\hfill \textbf{(b)} $k = 10$ \hfill\,\\
\includegraphics[width=0.49\textwidth]{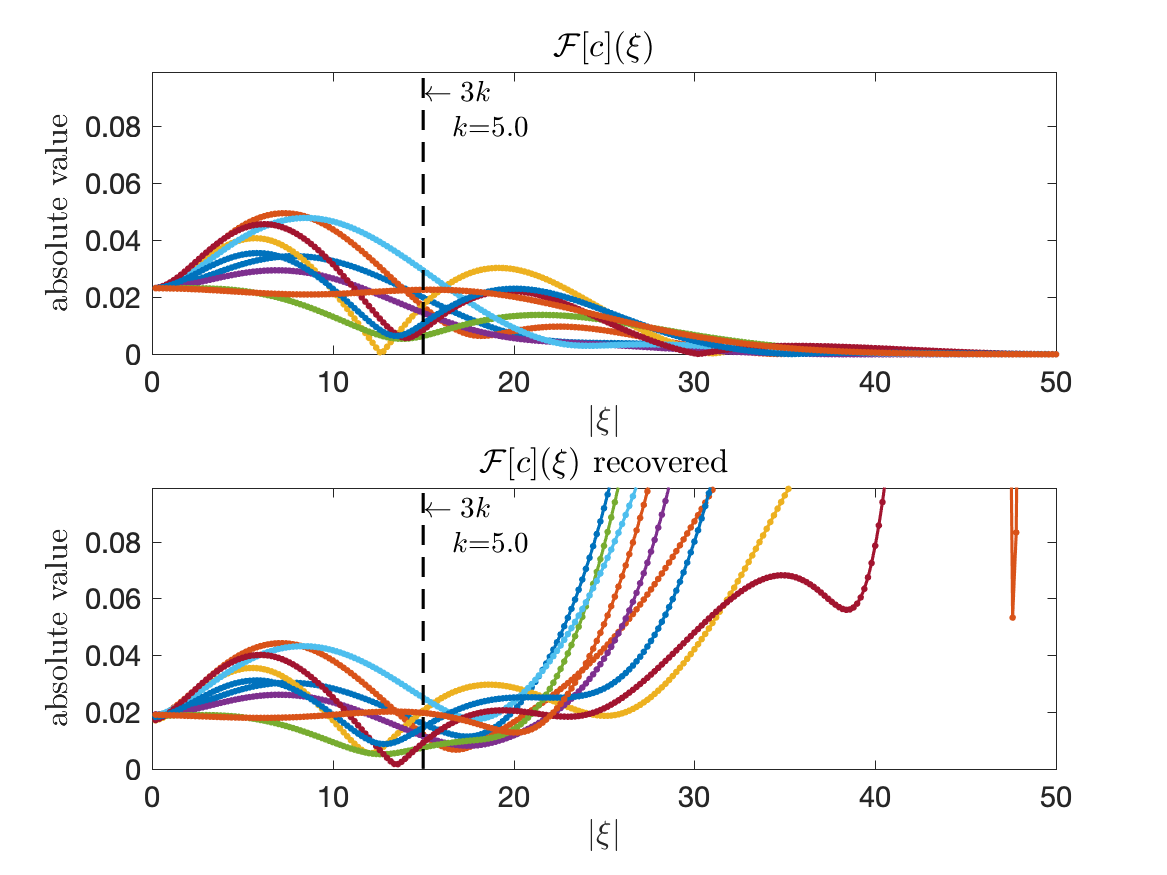}
\includegraphics[width=0.49\textwidth]{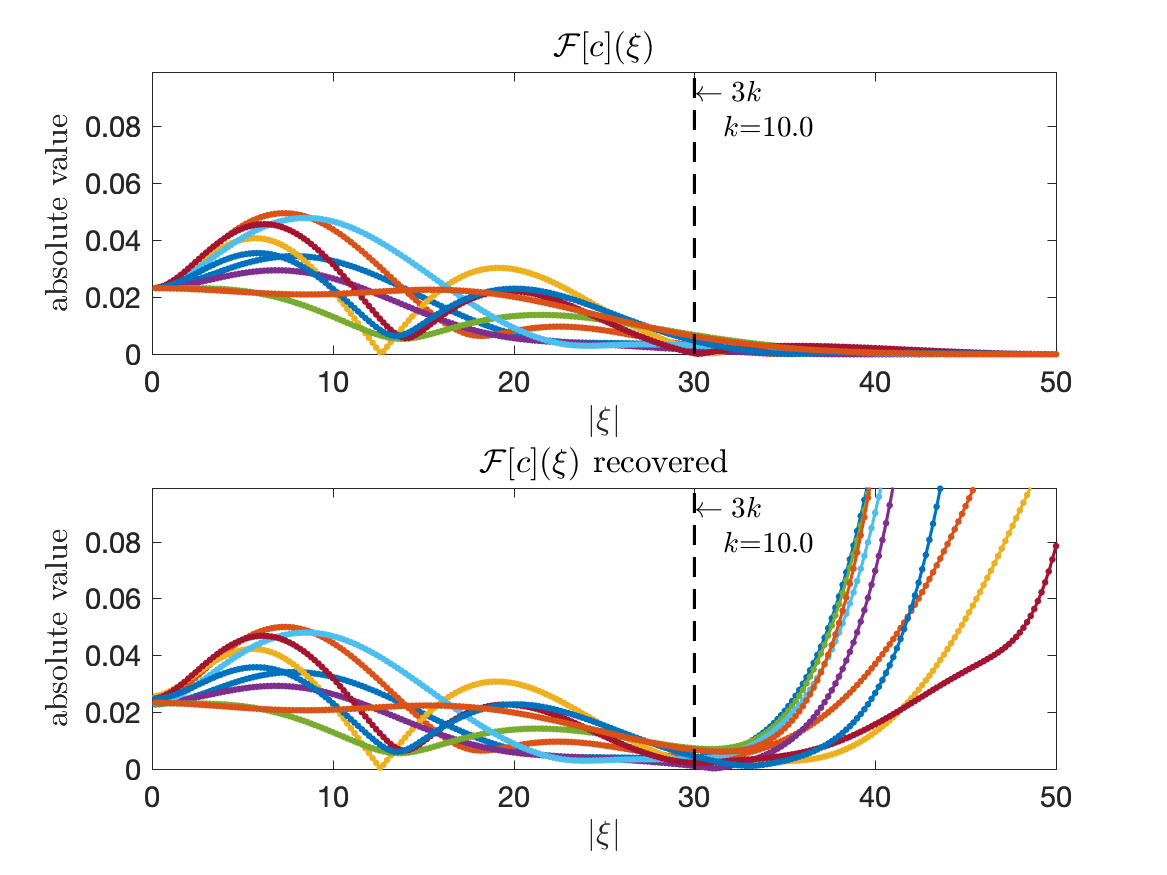}\\
\caption{The exact (Top) and recovered (Bottom) Fourier coefficients $\mathcal{F}[c](\xi)$ in each sub-figure: (a) $k=5$ and (b) $k=10$. Here, the horizontal axis shows the length $|\xi|$ of $\xi$; the vertical axis shows the absolute value $|\mathcal{F}[c](\xi)|$ of Fourier coefficients.}
\label{fig:Fc_inv}
\end{figure}

To numerically test \textbf{Algorithm 1}, we consider the domain $\Omega = B_{0.5}(0)$ in a square $[0.5,0.5]^{2}$. To avoid the inverse crime, we use a fine grids ($200 \times 200$ equal-distance points) for the forward problem and a coarse grid ($90 \times 90$ equal-distance points) for the inversion. The sampling points $\xi = (\xi_{1},\xi_{2})$ in frequency domain are shown in \Fref{fig:Xi}, marked by blue ``$\ast$'' near which all the Fourier coefficients will be recovered. In \Fref{fig:Fc_inv}, the horizontal axis shows the length $|\xi|$ of all $\xi$, and the vertical axis shows the absolute value $|\mathcal{F}[c](\xi)|$ of Fourier coefficients near the sampling points. By comparing the exact (top) and recovered (bottom) Fourier coefficients in each sub-figure of \Fref{fig:Fc_inv}: (a) $k = 5$ and (b) $k = 10$, we conclude that, while $k$ is larger, the more Fourier modes can be recovered stably, i.e.\ $\mathcal{F}[c](\xi)$ with $|\xi| \leq 3k$.

\begin{figure}[htbp]
\centering
\,\hfill \textbf{Quadratic case}: \hfill\,\\
\,\hfill \textbf{(a)} $k = 5$ \hfill\,\\
\includegraphics[width=0.6\textwidth]{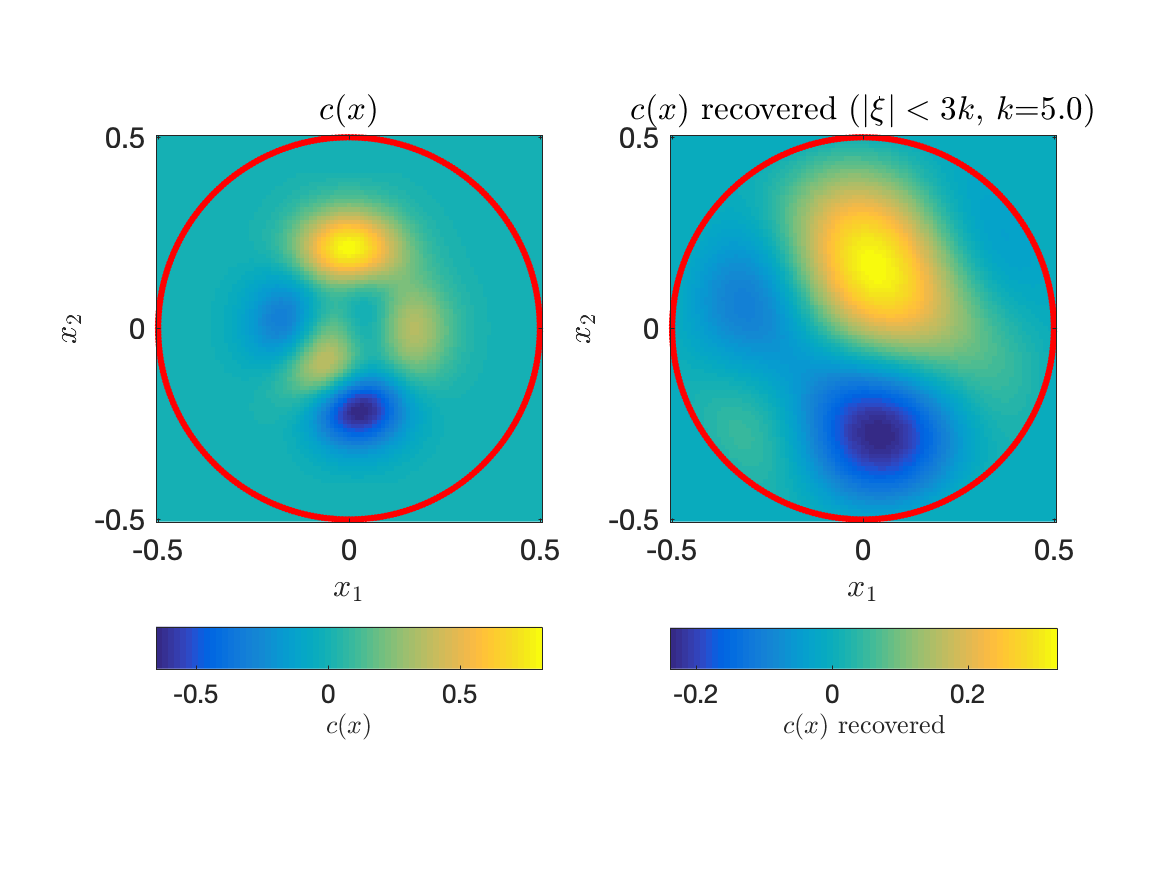}
\includegraphics[clip,trim={0.3in 0 3.6in 0},width=0.3\textwidth]{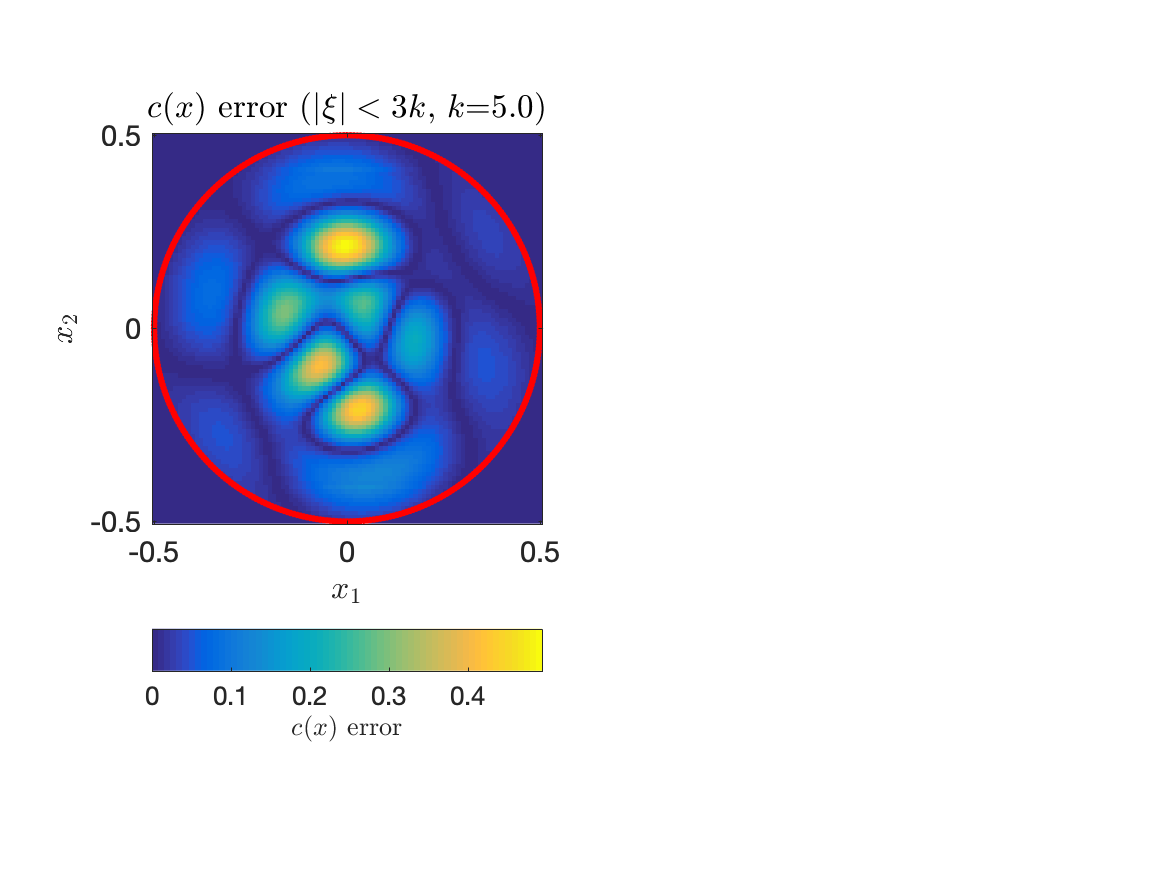}\\
\,\hfill \textbf{(b)} $k = 10$ \hfill\,\\
\includegraphics[width=0.6\textwidth]{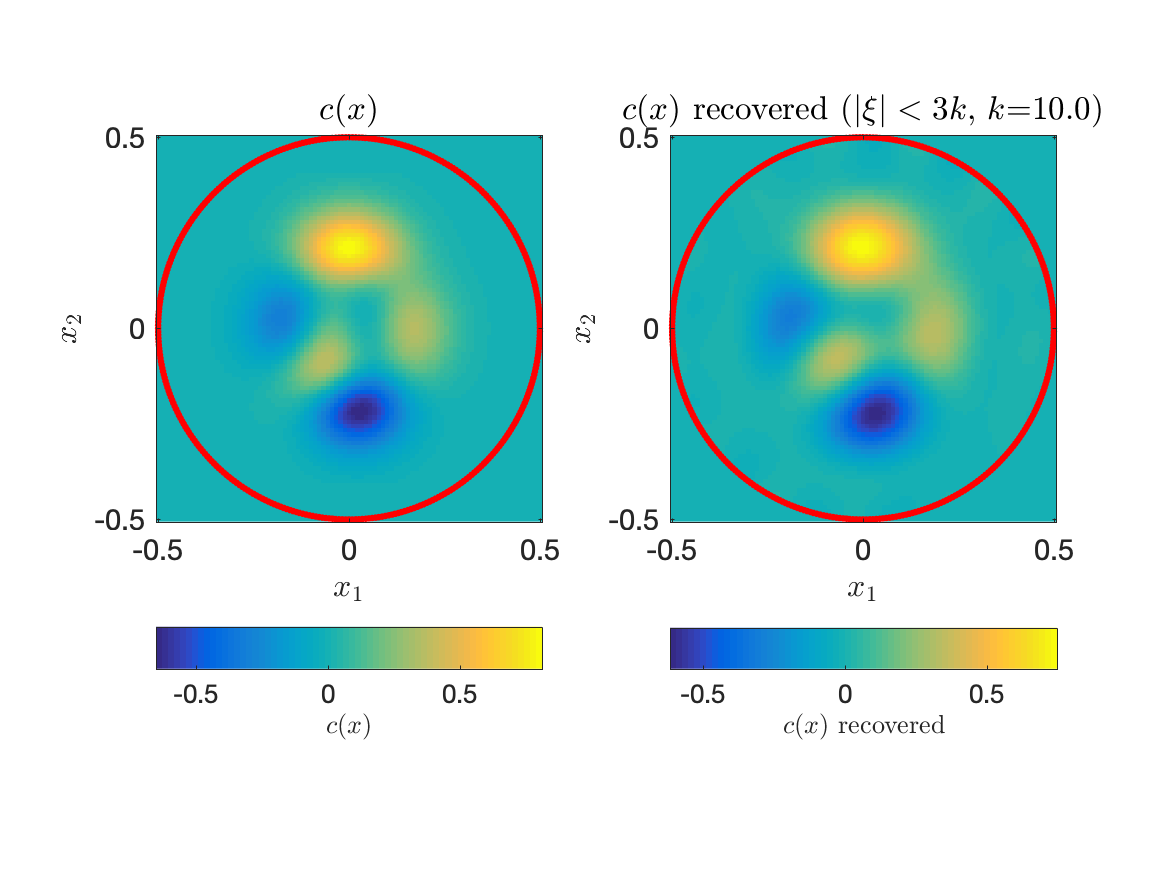}
\includegraphics[clip,trim={0.3in 0 3.6in 0},width=0.3\textwidth]{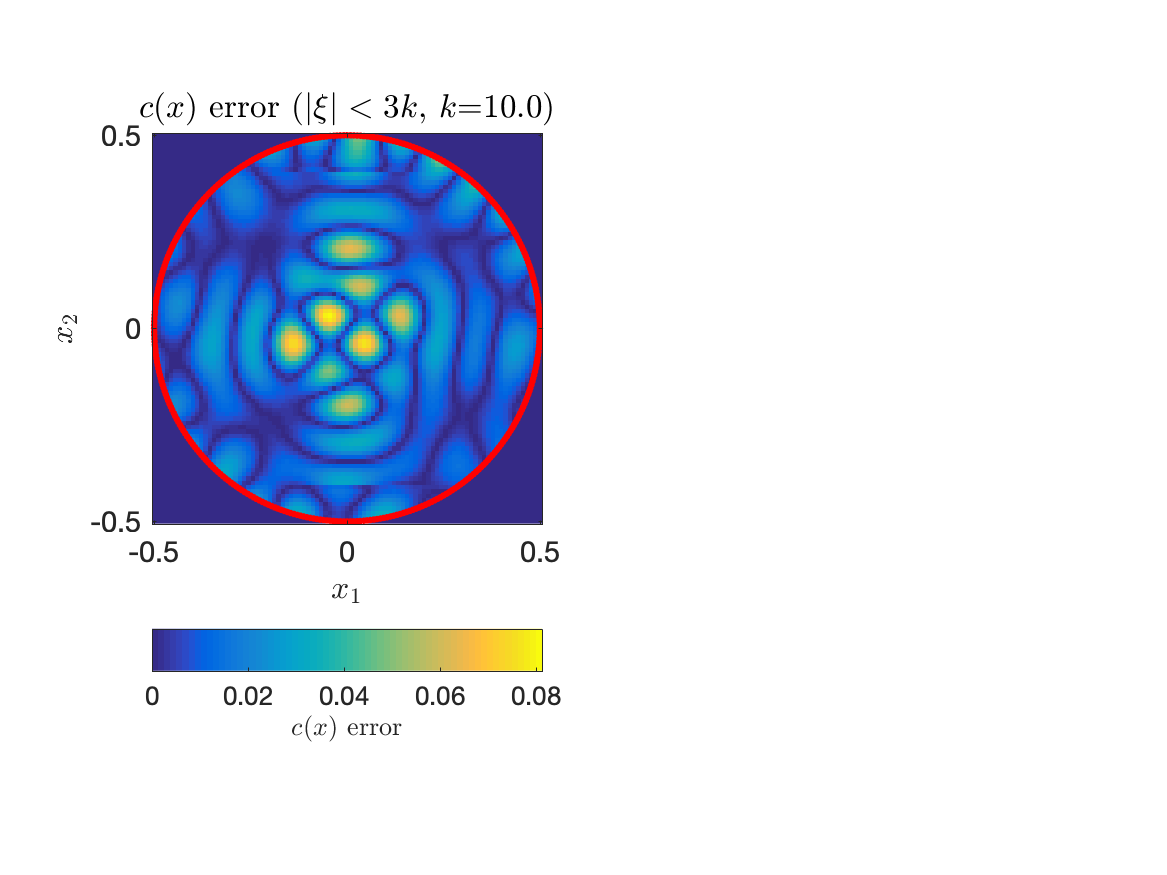}\\
\caption{\reviewA{The exact (Left) and recovered (Middle) potential functions $c(x)$ together with the point-wise absolute error (Right) when (a) $k=5$ and (b) $k=10$. Here, we use the Fourier coefficients $\mathcal{F}[c](\xi)$ in the range $|\xi| \leq 3k$.}}
\label{fig:Ic_inv}
\end{figure}

Then, by using all the recovered Fourier coefficients $\mathcal{F}[c](\xi)$ with $|\xi| \leq 3k$, we implement the inverse Fourier transform in Step 9 to reconstruct the potential function $c(x)$. In \Fref{fig:Ic_inv}, we present the exact and reconstructed potential functions $c(x)$ with different wavenumbers: (a) $k = 5$ and (b) $k = 10$, respectively. \reviewA{These results numerically verify the increasing stability in Theorem \ref{thm_quadratic} while $k$ becomes large. As one can observe, the point-wise absolute errors between the exact (left) and recovered (middle) potential functions are shown in \Fref{fig:Ic_inv}, and \textbf{Algorithm 1} reduces the maximum absolute error from $0.5$ to $0.08$ when $k$ increases from $5$ to $10$.}

\clearpage

\subsection{Reconstruction algorithm by $\Lambda'_{c}$ for high-order nonlinearity terms with multiple wavenumbers}\label{se4.2}

In this subsection, we show that the uniqueness result Theorem \ref{thm_unique_m} in \Sref{se2} indeed could provide a stable reconstruction algorithm for the nonlinear inverse Schr\"{o}dinger potential problem whose nonlinearity index is an arbitrary finite integer $m \geq 2$, if the linearized DtN map $\Lambda'_{c}$ of multiple wavenumbers is provided.

As highlighted in Remark \ref{rem_uniqueComSol}, the complex exponential solutions constructed in the proof of Theorem \ref{thm_unique_m} has a stable interval $\big[(m-1)k,(m+1)k\big]$ for any fixed $k$. Suppose that the same discrete (phase space) length and angle sets of the vectors in \Sref{se4.1} could be used, i.e.\ $\{\kappa_{i}\}_{i=1}^{I}$, $\{\hat{y}_{s}\}_{s=1}^{S}$, $\{\hat{z}_{s}\}_{s=1}^{S}$ and the following vectors $\mu_{\ell}^{\langle i;s \rangle} \in \mathbb{C}^{n}$, $\ell = 1,2$ are chosen
\begin{equation*}
\left\{~
\eqalign{
\mu_{1}^{\langle i;s \rangle} &:= \frac{+(m^{2}-1)k^{2}+\kappa_{i}^{2}}{2m\kappa_{i}} \hat{y}_{s} - \frac{\sqrt{-(m^{2}-1)^{2}k^{4}+2(m^{2}+1)k^{2}\kappa_{i}^{2}-\kappa_{i}^{4}}}{2m\kappa_{i}} \hat{z}_{s}, \\
\mu_{2}^{\langle i;s \rangle} &:= \frac{-(m^{2}-1)k^{2}+\kappa_{i}^{2}}{2\kappa_{i}} \hat{y}_{s} + \frac{\sqrt{-(m^{2}-1)^{2}k^{4}+2(m^{2}+1)k^{2}\kappa_{i}^{2}-\kappa_{i}^{4}}}{2\kappa_{i}} \hat{z}_{s},
}
\right.
\end{equation*}
similar to \textbf{Algorithm 1}, we summarize a plain reconstruction algorithm for high-order nonlinearity terms with a fixed wavenumber $k$, according to the identity \eref{eq_calderonILX}.

\vspace{10pt}
\hrule\hrule
\vspace{8pt}
{\parindent 0pt \bf Algorithm 2: Reconstruction Algorithm for the Linearized Schr\"{o}dinger Potential Problem, the high-order nonlinearity term} %
\vspace{5pt}
\hrule
\vspace{8pt}
{\parindent 0pt \bf Input:} %
$m$, $k$, %
$\{\kappa_{i}\}_{i=1}^{I}$, %
$\{\hat{y}_{s}\}_{s=1}^{S}$, $\{\hat{z}_{s}\}_{s=1}^{S}$ and %
$\{\sigma^{\langle i;s \rangle}\}$; \\[5pt]%
{\parindent 0pt \bf Output:} %
Approximated Potential $c^{\langle I+1;1 \rangle}$. \\[-15pt]%
\begin{enumerate}
  \item[1:] $\,$ Set $c^{\langle 1;1 \rangle} := 0$; %
  \item[2:] $\,$ {\bf For} $i = 1,2,\dots,I$ (length~updating) %
  \item[3:] $\,$ \quad {\bf For} $s = 1,2,\dots,S$ (angle~updating) %
  \item[4:] $\,$ \quad \quad Choose $u_{0} := \exp \{ \rmi \mu_{1}^{\langle i;s \rangle} \cdot x \}$; %
  \item[5:] $\,$ \quad \quad Measure the Neumann boundary data $\partial_{\nu} u$ of the forward problem \eref{eq_HelmholtzEqmain} %
  \item[] $\,$ \quad \qquad while the Dirichlet boundary data $u_{0}|_{\partial\Omega}$ are given; %
  \item[6:] $\,$ \quad \quad Calculate the approximated linearized Neumann boundary data %
  \item[] $\,$ \quad \qquad $g_{u}^{\,\prime} := (\partial_{\nu} u - \partial_{\nu} u_{0})|_{\partial\Omega}$; %
  \item[7:] $\,$ \quad \quad Choose $\varphi := \exp \{ \rmi \mu_{2}^{\langle i;s \rangle} \cdot x \}$ and $\gamma := \big[ u_{0}^{m} \varphi \big]^{-1} = \exp \{ -\rmi \xi^{\langle i;s \rangle} \cdot x \}$; %
  \item[8:] $\,$ \quad \quad Compute $\mathcal{F}[c](\xi^{\langle i;s \rangle}) \approx \int_{\partial\Omega} g_{u}^{\,\prime} \, \varphi \,\rmd S$; %
  \item[9:] $\,$ \quad \quad Update $c^{\langle i;s+1 \rangle} := c^{\langle i;s \rangle} + \mathcal{F}[c](\xi^{\langle i;s \rangle}) \, \gamma \sigma^{\langle i;s \rangle}$, \quad if $\kappa_{i} \in \big[(m-1)k,(m+1)k\big]$; %
  \item[10:] $\,$ \quad {\bf End}; %
  \item[11:] $\,$ \quad Set $c^{\langle i+1;1 \rangle} := c^{\langle i;S+1 \rangle}$; %
  \item[12:] $\,$ {\bf End}. %
\end{enumerate}
\vspace{8pt}
\hrule\hrule
\vspace{10pt}

Furthermore, if we could measure the boundary data by appropriate multiple wavenumbers, we could reconstruct sufficiently many Fourier coefficients of the unknown potential function. By choosing $k_{1}$ small and a threshold value $K$ as the maximum wavenumber, we choose a discrete set of multiple wavenumbers, namely
\begin{equation}\label{eq_numerMultiwavenumber}
\{ k_{j} \}_{j=1}^{J} \subset (0, K \,],
\end{equation}
which satisfies $k_{j+1} = \frac{m+1}{m-1} k_{j}$. Below we present an updated reconstruction algorithm of \textbf{Algorithm 2} for the linearized Schr\"{o}dinger potential problem with a high-order nonlinearity term, i.e.\ the nonlinearity index $m \geq 2$, if the linearized DtN map $\Lambda'_{c}$ of multiple wavenumbers can be obtained.

\vspace{10pt}
\hrule\hrule
\vspace{8pt}
{\parindent 0pt \bf Algorithm 2*: Reconstruction Algorithm for the Linearized Schr\"{o}dinger Potential Problem with a high-order nonlinearity term (Multiple wavenumbers)} %
\vspace{5pt}
\hrule
\vspace{8pt}
{\parindent 0pt \bf Input:} %
$m$, %
$\{ k_{j} \}_{j=1}^{J}$, %
$\{\kappa_{i}\}_{i=1}^{I}$, %
$\{\hat{y}_{s}\}_{s=1}^{S}$, $\{\hat{z}_{s}\}_{s=1}^{S}$ and %
$\{\sigma^{\langle i;s \rangle}\}$; \\[5pt]%
{\parindent 0pt \bf Output:} %
Approximated Potential $c_{\rm inv} := \sum\limits_{j=1}^{J} c^{\langle I+1;1 \rangle}_{j}$. \\[-15pt]%
\begin{enumerate}
  \item[1:] $\,$ {\bf For} $j = 1,2,\dots,J$ (wavenumber~updating) %
  \item[2:] $\,$ \quad Compute the approximated potential $c^{\langle I+1;1 \rangle}_{j}$ by using \textbf{Algorithm 2} and a fixed $k_{j}$; %
  \item[3:] $\,$ {\bf End}. %
\end{enumerate}
\vspace{8pt}
\hrule\hrule
\vspace{10pt}

\begin{figure}[htbp]
\centering
\,\hfill \textbf{Cubic case}: \hfill\,\\
\,\hfill \textbf{multiple wavenumbers} $k \in \{ 1.25, 2.5, 5, 10 \}$ \hfill\,\\
\includegraphics[width=0.6\textwidth]{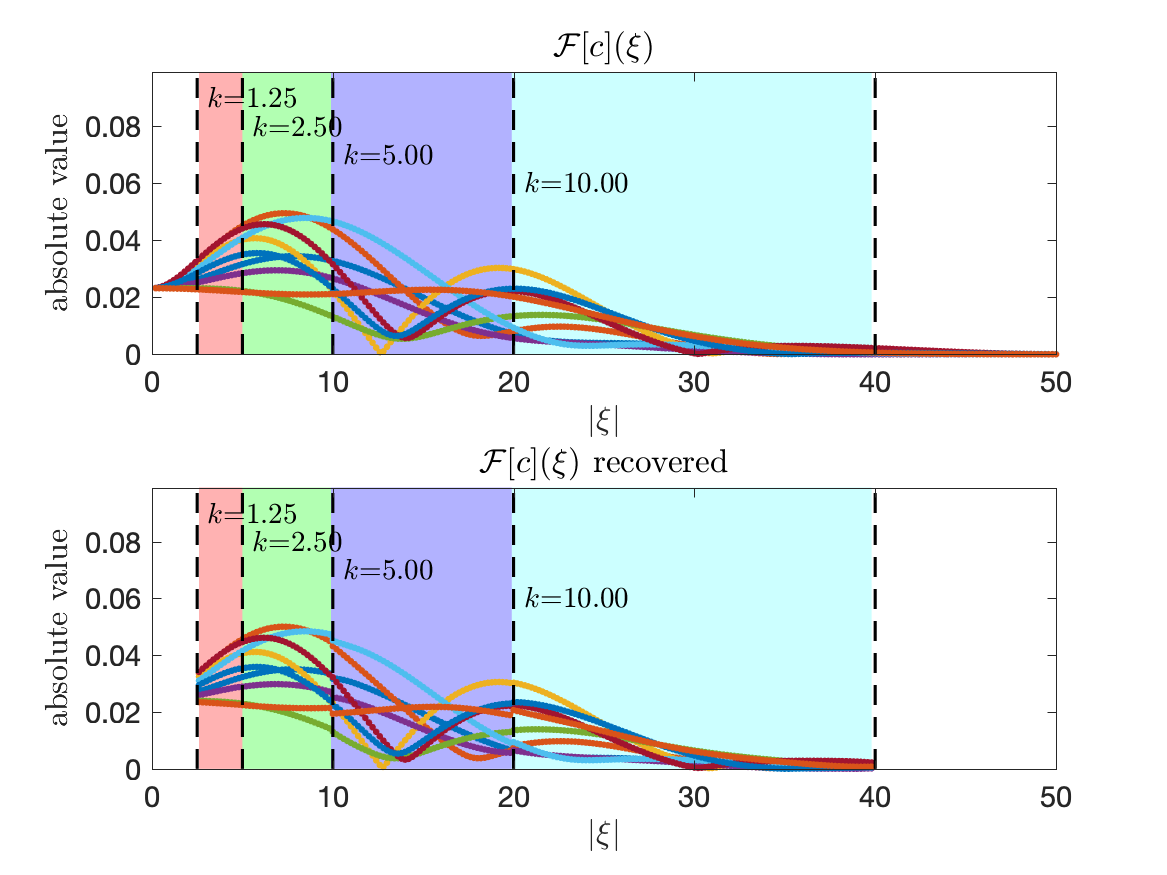}\\
\caption{(Cubic case, $m = 3$) The exact (Top) and recovered (Bottom) Fourier coefficients $\mathcal{F}[c](\xi)$ with multiple wavenumbers $k \in \{ 1.25, 2.5, 5, 10 \}$. Here, the horizontal axis shows the length $|\xi|$ of $\xi$; the vertical axis shows the absolute value $|\mathcal{F}[c](\xi)|$ of Fourier coefficients.}
\label{fig:3_Fc_inv_cmb}
\end{figure}

As an illustration, we consider the linearized Schr\"{o}dinger potential problem with a cubic nonlinear term ($m = 3$). The wavenumber set in \eref{eq_numerMultiwavenumber} is set with $k_{1} = 1.25$ and $K = 10$ where we recover the Fourier coefficients $\mathcal{F}[c](\xi)$ with $4$ wavenumbers $k \in \{ 1.25, 2.5, 5, 10 \}$. In \Fref{fig:3_Fc_inv_cmb}, the red region indicates the Fourier coefficients within $[2k_{1},4k_{1}) = [2.5,5)$, the green region indicates the Fourier modes within $[2k_{2},4k_{2}) = [5,10)$, the blue region indicates the Fourier modes within $[2k_{3},4k_{3}) = [10,20)$, and the cyan region indicates the Fourier coefficients within $[2k_{4},4k_{4}) = [20,40)$.

\begin{figure}[htbp]
\centering
\,\hfill \textbf{Cubic case}: \hfill\,\\
\,\hfill \textbf{multiple wavenumbers} $k \in \{ 1.25, 2.5, 5, 10 \}$ \hfill\,\\
\includegraphics[width=0.6\textwidth]{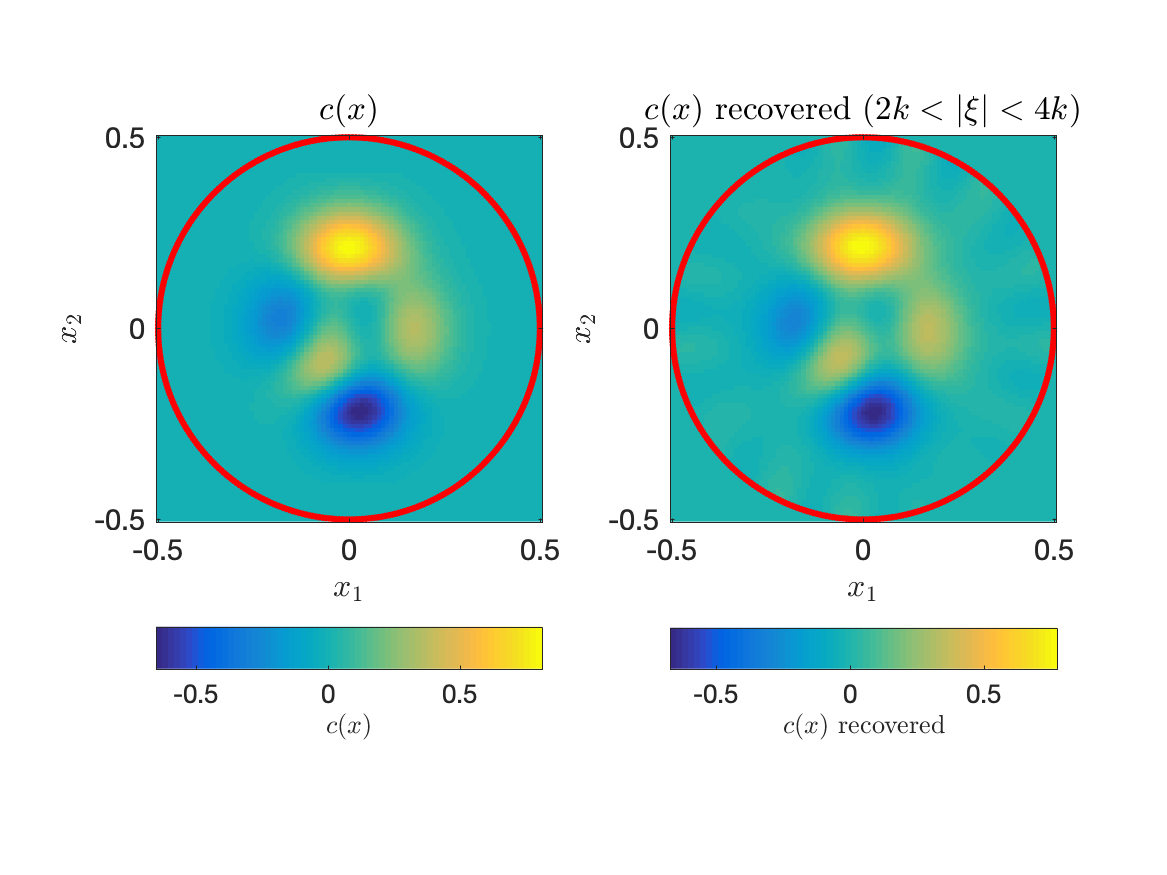}
\includegraphics[clip,trim={0.3in 0 3.6in 0},width=0.3\textwidth]{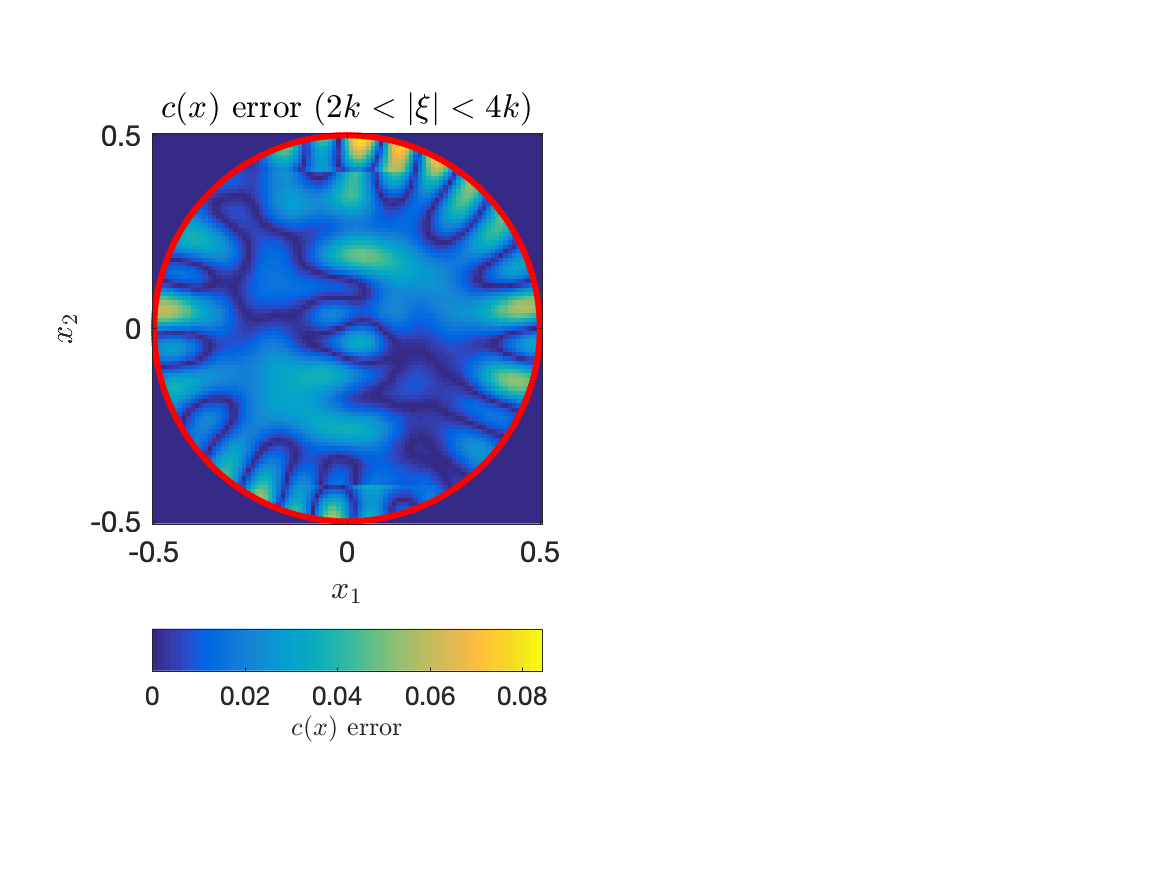}\\
\caption{(Cubic case, $m = 3$) \reviewA{The exact (Left) and recovered (Middle) potential functions $c(x)$ together with the point-wise absolute error (Right) when multiple wavenumbers $k \in \{ 1.25, 2.5, 5, 10 \}$ are considered.}}
\label{fig:3_Ic_inv_cmb}
\end{figure}

By using Fourier coefficients $\mathcal{F}[c](\xi)$ within $|\xi| \in \bigcup_{j=1}^{J} \big[(m-1)k_{j},(m+1)k_{j}\big) = [(m-1)k_{1},(m+1)k_{J})$, we implement the inverse Fourier transform to reconstruct the potential function $c(x)$. In \Fref{fig:3_Ic_inv_cmb}, we present the exact (left) and reconstructed (right) potential functions $c(x)$ with $4$ wavenumbers $k \in \{ 1.25, 2.5, 5, 10 \}$. It can be seen that, by including the boundary measurements of four wavenumbers, we have obtained a good approximation of the unknown potential function in \eref{eq_HelmholtzEqmain} with a cubic nonlinear term $m = 3$.

\subsection{Vanilla reconstruction algorithm by $D^{m}_{0}\Lambda_{c}$}\label{se4.3}

Noticing that the linearized DtN map $\Lambda'_{c}$ can be approximated by ignoring the high order terms, i.e.\ \eref{eq_se4linearizedNeumann}, we are allowed to adopt this idea and design a vanilla reconstruction algorithm for another linearized DtN map $D^{m}_{0}\Lambda_{c}$.

As shown in \Sref{se1.2}, the identity \eref{eq_nonlinearsmallboundary_equality} plays a key role in recovering the unknown potential function $c(x)$ with respect to the small boundary data $f_{\varepsilon}$. More precisely, it relies on the boundary data $\partial_{\nu} w |_{\partial\Omega} = D^{m}_{0} \Lambda_{c}(f_{1},\ldots,f_{m})$ sensitively. Thus in this subsection, we provide some numerical tests to study the consequence by the linearized DtN map $D^{m}_{0} \Lambda_{c}$, and the following formulae are employed to approximate the $m$th Fr\'{e}chet derivative $D^{m}$ with $m = 2$ and $3$, such that
\begin{equation}\label{eq_se4appro}
\left\{~
\eqalign{
\partial_{\varepsilon_{1}} \partial_{\varepsilon_{2}} \Lambda_{c} (f_{\varepsilon}) &\approx \frac{1}{\varepsilon_{1}\varepsilon_{2}} \Big( \Lambda_{c} (\varepsilon_{1} f_{1}+\varepsilon_{2} f_{2}) - \Lambda_{c} (\varepsilon_{2} f_{2}) - \Lambda_{c} (\varepsilon_{1} f_{1}) + \Lambda_{c} (0) \Big), \\
\partial_{\varepsilon_{1}} \partial_{\varepsilon_{2}} \partial_{\varepsilon_{3}} \Lambda_{c} (f_{\varepsilon}) &\approx \frac{1}{\varepsilon_{1}\varepsilon_{2}\varepsilon_{3}} \Big(
\Lambda_{c} (\varepsilon_{1} f_{1}+\varepsilon_{2} f_{2}+\varepsilon_{3} f_{3}) \\
&\quad - \Lambda_{c} (\varepsilon_{1} f_{1}+\varepsilon_{2} f_{2})
- \Lambda_{c} (\varepsilon_{1} f_{1}+\varepsilon_{3} f_{3})
- \Lambda_{c} (\varepsilon_{2} f_{2}+\varepsilon_{3} f_{3}) \\
&\quad + \Lambda_{c} (\varepsilon_{3} f_{3})
+ \Lambda_{c} (\varepsilon_{2} f_{2})
+ \Lambda_{c} (\varepsilon_{1} f_{1})
- \Lambda_{c} (0) \Big),
}
\right.
\end{equation}
when each $\varepsilon_{j}$, $j=1,2,3$ is small enough and chosen appropriately. Here we mention that $\Lambda_{c} (0) = 0$.

\begin{figure}[htbp]
\centering
\,\hfill \textbf{Linearized DtN map $D^{2}_{0} \Lambda_{c}$ (Quadratic case)}: \hfill\,\\
\includegraphics[width=0.45\textwidth]{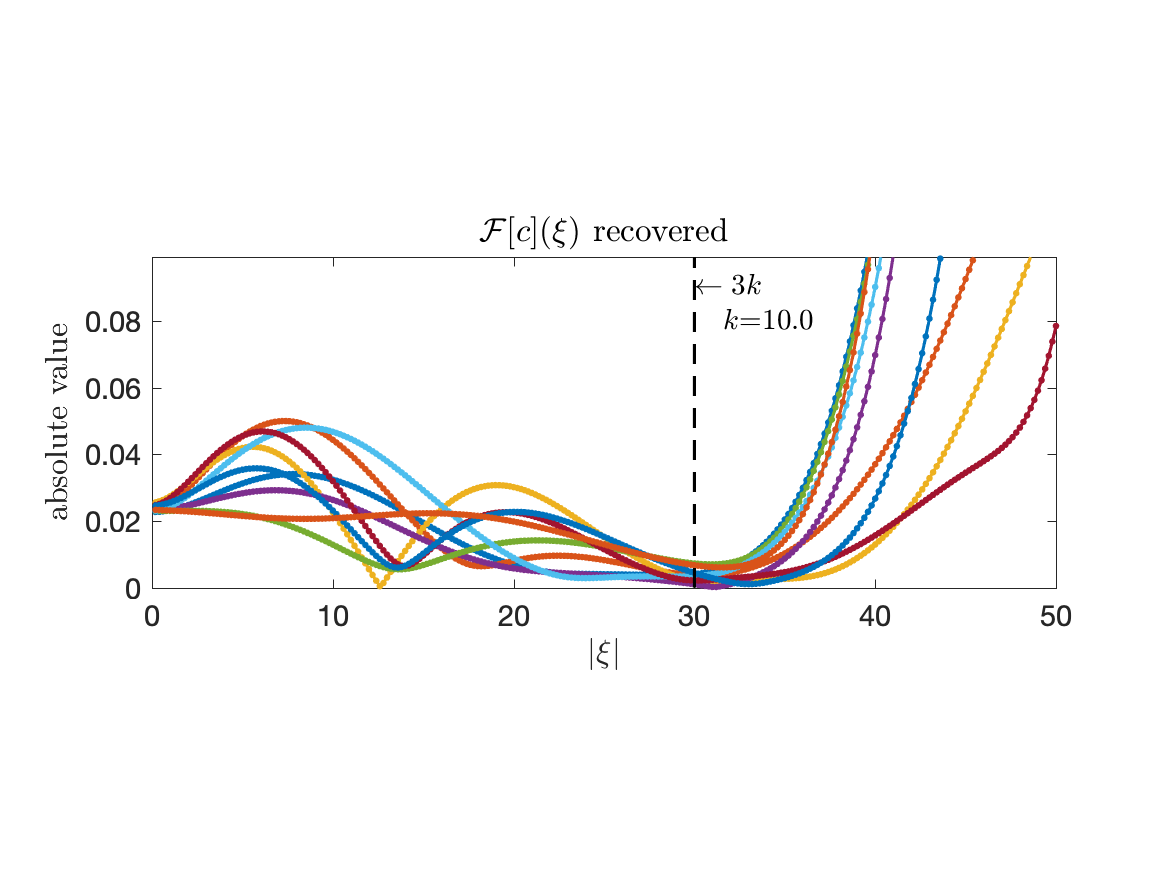}
\includegraphics[width=0.45\textwidth]{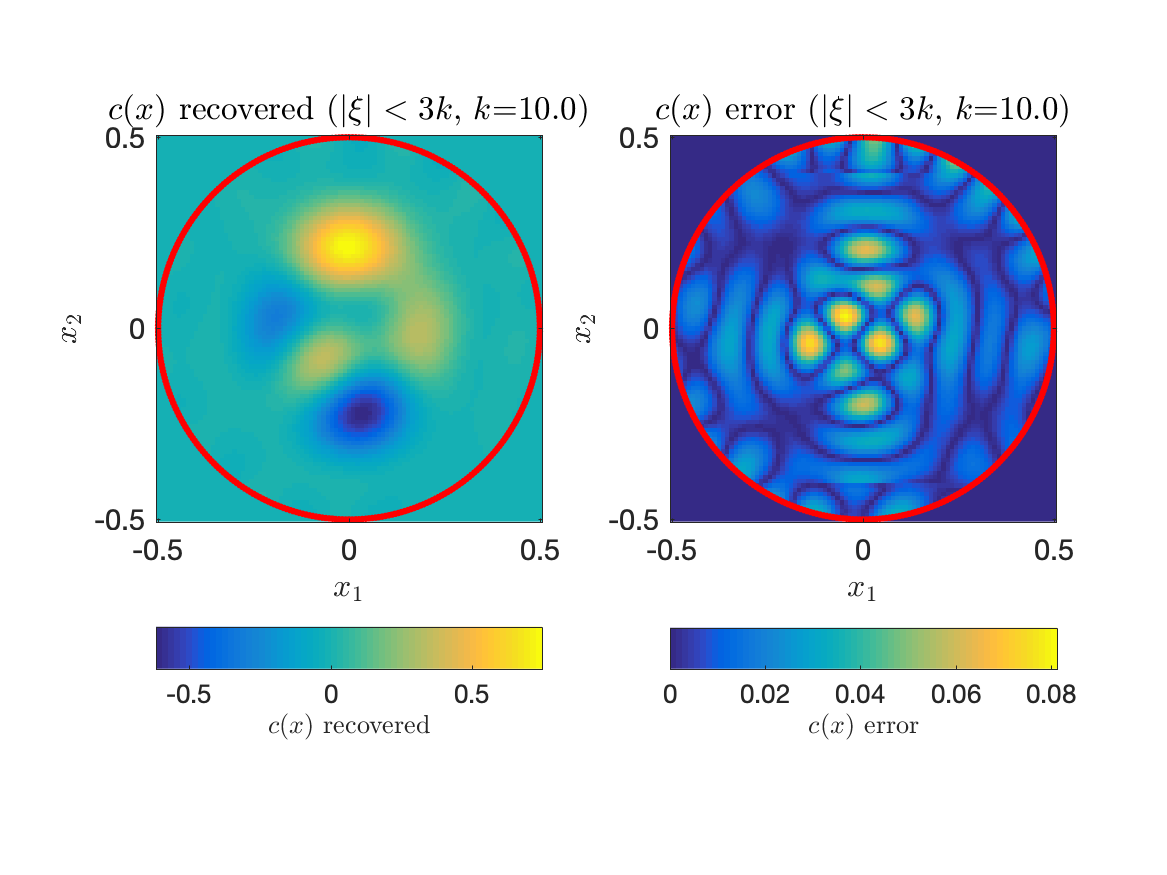}\\
\,\hfill \textbf{Linearized DtN map $D^{3}_{0} \Lambda_{c}$ (Cubic case)}: \hfill\,\\
\includegraphics[width=0.45\textwidth]{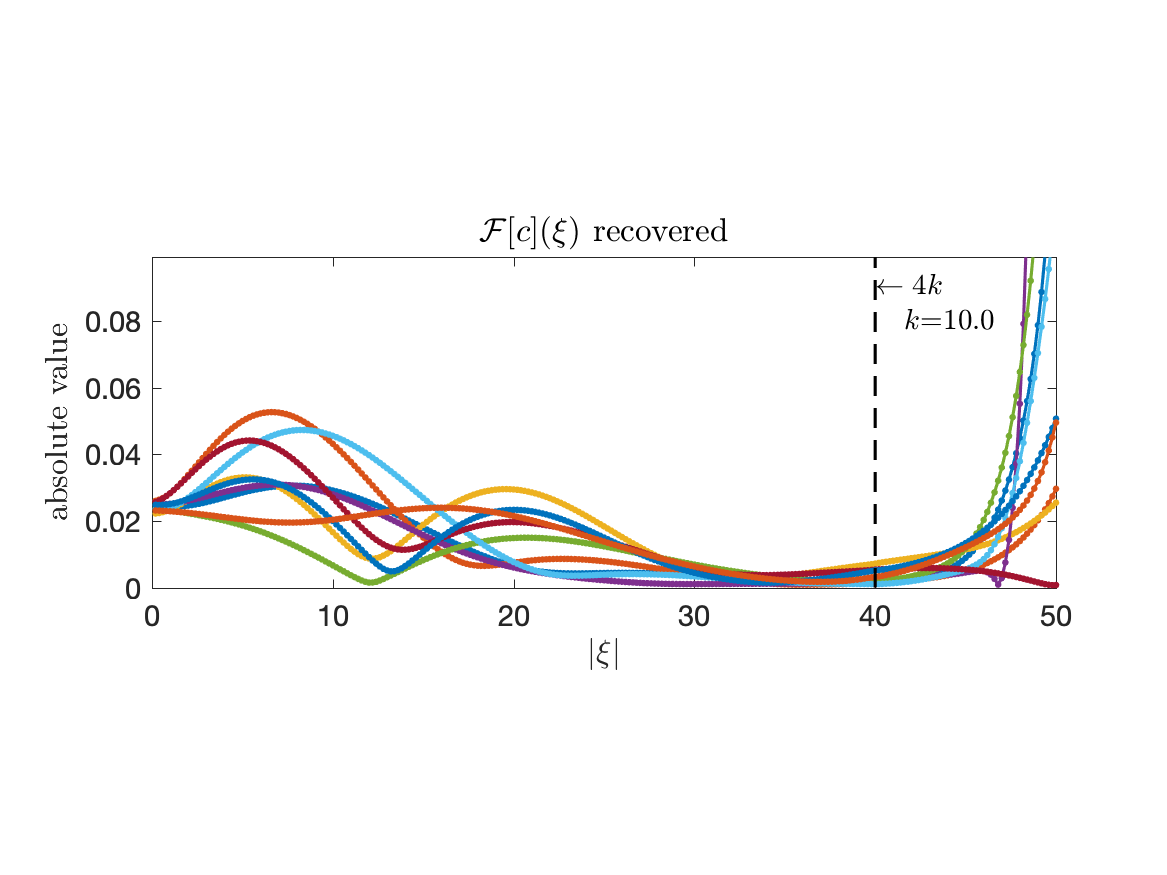}
\includegraphics[width=0.45\textwidth]{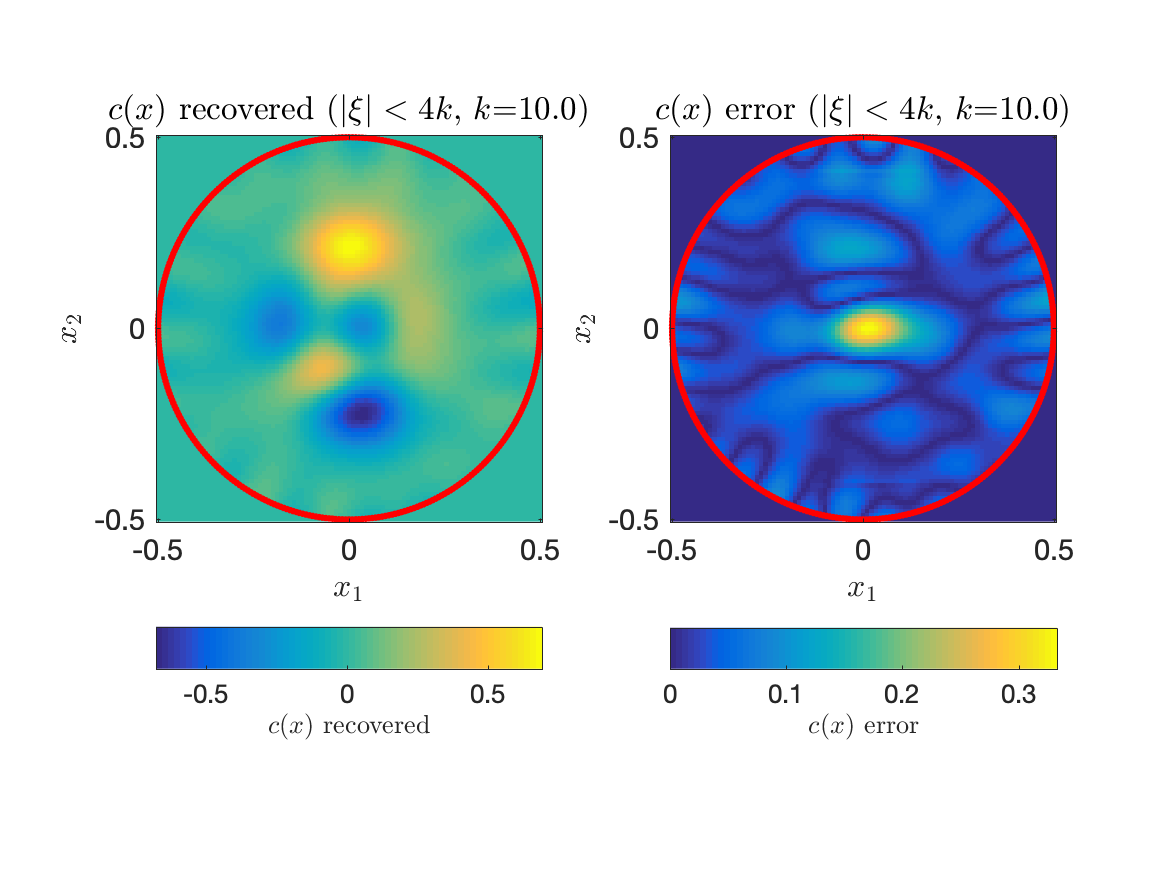}\\
\caption{Left: The recovered Fourier coefficients $\mathcal{F}[c](\xi)$ with $k = 10$ by linearized DtN map $D^{m}_{0} \Lambda_{c}$. \reviewA{Middle: The recovered potential $c(x)$ with $k = 10$ and $|\xi| \leq (m+1)k$ by linearized DtN map $D^{m}_{0} \Lambda_{c}$. Right: the point-wise absolute error between the exact and recovered potential functions.} Here $m=2$ (Top) and $m=3$ (Bottom).}
\label{fig:MS_Ic_inv}
\end{figure}

We note that one can modify {\bf Algorithm 1} carefully to design a reconstruction algorithm for the linearized DtN map $D^{m}_{0}\Lambda_{c}$ if appropriate complex exponential solutions \eref{eq_thm1_CExsol1} or \eref{eq_thm1_CExsol2} in the proof of Theorem \ref{thm_holder} are chosen and the above derivative approximation schemes \eref{eq_se4appro} are implemented. To save the space, we skip the pseudocode of the algorithm but present the reconstructed potential function $c(x)$ and its Fourier coefficients $\mathcal{F}[c](\xi)$ in \Fref{fig:MS_Ic_inv} for different nonlinearity index with $m = 2,3$. In both cases, we have chosen $\varepsilon_{j} = 0.1$, $j=1,2,3$ as illustration. In principle, one can extend the derivative approximation formulae \eref{eq_se4appro} to more general case with $m > 3$ and tune the small parameters $\varepsilon_{j}$ carefully to obtain better resolution. But this is beyond the scope of current work and will be considered as future work.

\subsection{\reviewC{The noise propagation}}

\begin{figure}[htbp]
\centering
\,\hfill \textbf{Quadratic case with noise (linearized DtN map $\Lambda'_{c}$)}: \hfill\,\\
\includegraphics[width=0.45\textwidth]{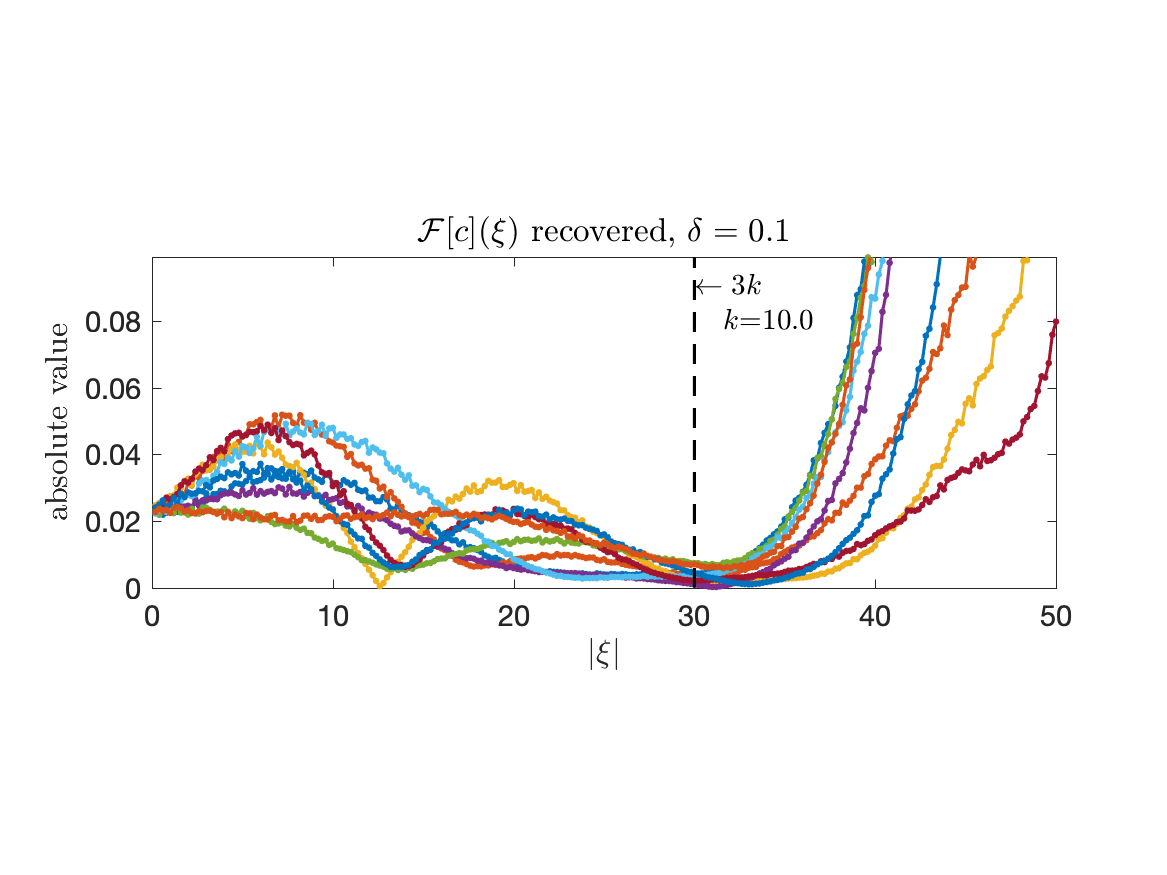}
\includegraphics[width=0.45\textwidth]{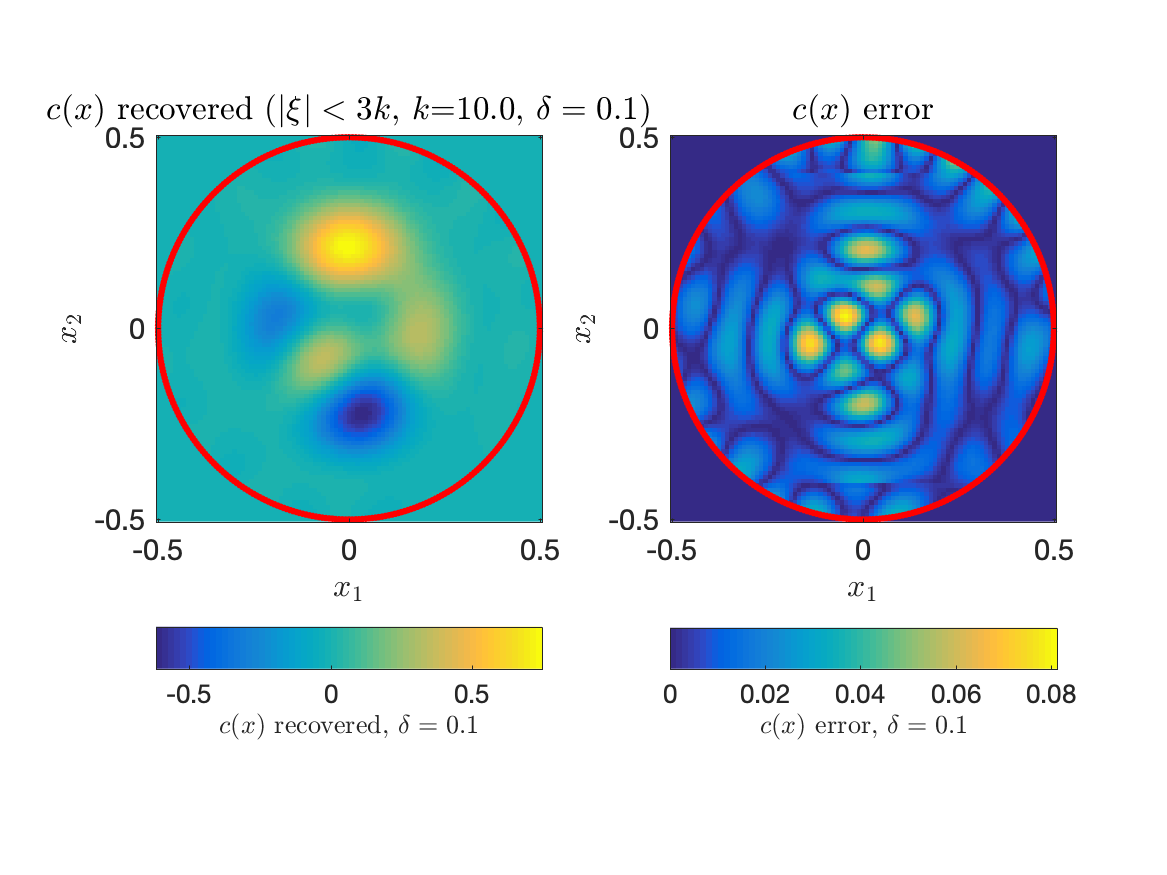}\\
\,\hfill \textbf{Cubic case with noise (linearized DtN map $D^{m}_{0} \Lambda_{c}$)}: \hfill\,\\
\includegraphics[width=0.45\textwidth]{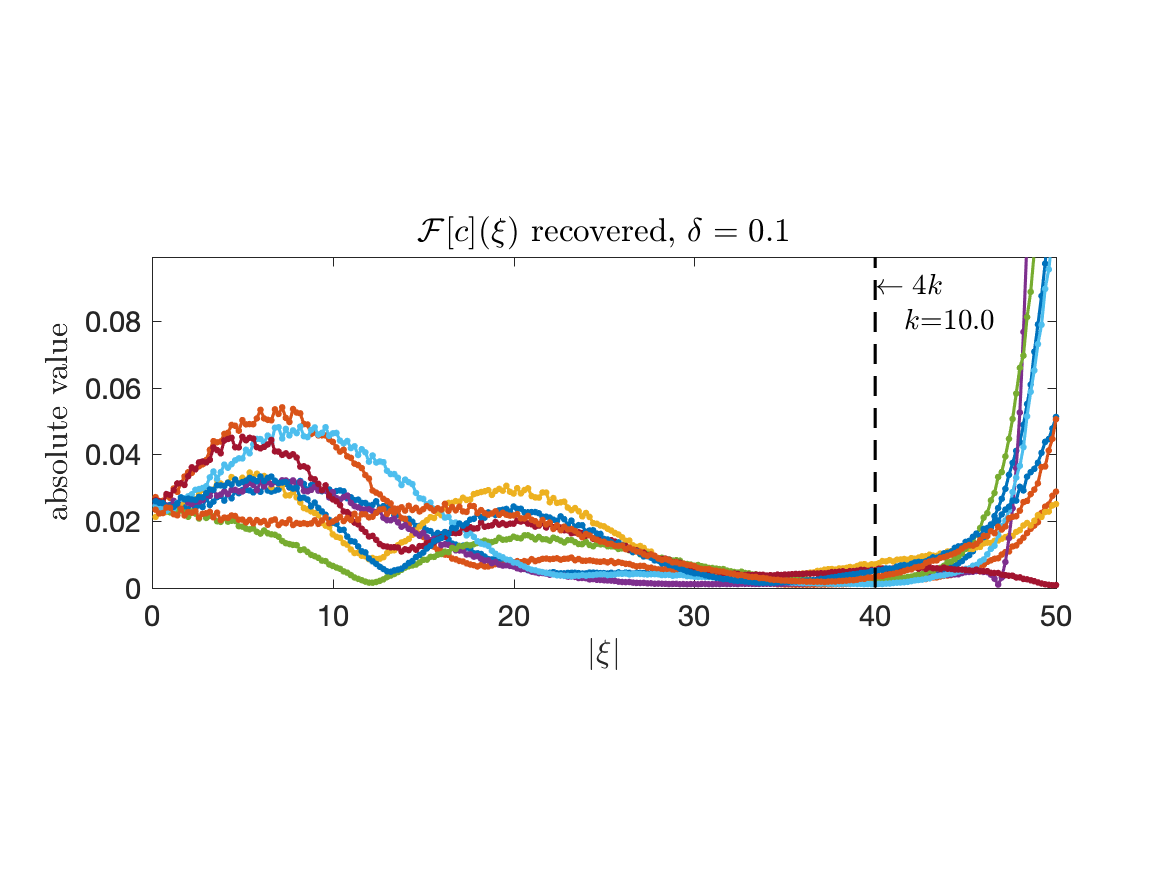}
\includegraphics[width=0.45\textwidth]{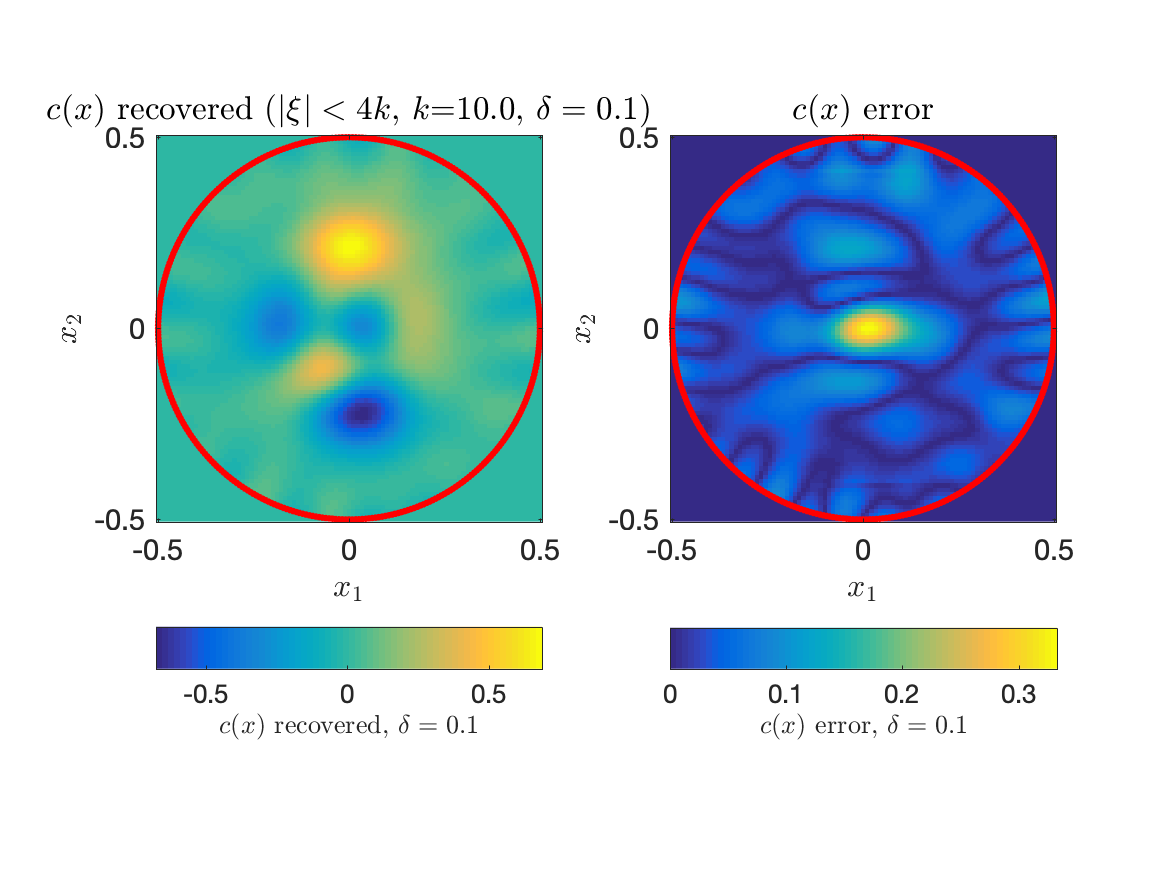}\\
\caption{\reviewC{Left: The recovered Fourier coefficients $\mathcal{F}[c](\xi)$ with $k = 10$ and $\delta = 0.1$. Middle: The recovered potential $c(x)$ with $k = 10$, $|\xi| \leq (m+1)k$ and $\delta = 0.1$. Right: The point-wise absolute error between the exact and recovered potential functions. Here, the quadratic case $m=2$ by linearized form $\Lambda'_{c}$ (Top) and the cubic case $m=3$ by linearized form $D^{m}_{0} \Lambda_{c}$ (Bottom).}}
\label{fig:noisy_Ic_inv}
\end{figure}

\reviewC{Finally, we consider the noise propagation on both linearized Neumann boundary data $\partial_{\nu} w |_{\partial\Omega} = D^{m}_{0} \Lambda_{c}(f_{1},\ldots,f_{m})$ in \eref{eq_HighOrderLinearDtN} and $\partial_{\nu} u_{1}|_{\partial\Omega} = \Lambda'_{c} g_{0}$ in \eref{eq_LinearDtN}.
Assume that there exists a (relative) noise level $\delta$ such that the uniformly bounded noise between the exact and noisy linearized Neumann boundary data satisfies
\begin{equation*}
\frac{\left\| (\partial_{\nu} w)^{\delta} - \partial_{\nu} w \right\|_{L^{\infty}(\partial\Omega)}}{\left\| \partial_{\nu} w \right\|_{L^{\infty}(\partial\Omega)}} \leqslant \delta, 
\qquad
\frac{\left\| (\partial_{\nu} u_{1})^{\delta} - \partial_{\nu} u_{1} \right\|_{L^{\infty}(\partial\Omega)}}{\left\| \partial_{\nu} u_{1} \right\|_{L^{\infty}(\partial\Omega)}} \leqslant \delta,
\end{equation*}
where $(\partial_{\nu} w)^{\delta}$ and $(\partial_{\nu} u_{1})^{\delta}$ denote the noisy Neumann boundary data, respectively.
We present the recovered Fourier coefficients (left) and potential function (middle) together with the point-wise absolute error (right) in each sub-figure of Figure \ref{fig:noisy_Ic_inv}, where $k = 10$ and $\delta = 0.1$. Though the recovered Fourier coefficients become rough when noise appears, the recovered potential retains good resolution in both linearized DtN maps $\Lambda'_{c}$ and $D^{m}_{0} \Lambda_{c}$.
More specifically, the results in \Fref{fig:noisy_Ic_inv}, recovered by the noisy measurements, can be compared with the corresponding noiseless results in \Fref{fig:Fc_inv}(b), \Fref{fig:Ic_inv}(b) for the quadratic case $m = 2$ by $\Lambda'_{c}$, and in \Fref{fig:MS_Ic_inv} (bottom) for the cubic case $m = 3$ by $D^{m}_{0} \Lambda_{c}$.
Indeed it can be observed that, with the chosen truncated value $(m+1)k$, the reconstructed results are robust with respect to the noise propagation because of our chosen complex exponential solutions in both cases.}

\clearpage




\begin{thebibliography}{99}

%
%
%
%
%
%
%
%
%
%
%
%
%
%
%
%
%
%
%
%
%
%
%
%
%
%
%
%
%
%
%
%
%
%
%
%
%
%

\bibitem{A1988}
Alessandrini, G {\it Stable determination of conductivity by boundary measurement}. Appl. Anal. {\bf 27} (1988), 153--172.


\bibitem{BG2015}
Borges, C; Greengard, L {\it Inverse obstacle scattering in two dimensions with multiple frequency data and multiple angles of incidence}. SIAM J. Imaging Sci. {\bf 8} (2015), no. 1, 280--298.


\bibitem{BLLT2015}
Bao, G; Li, P; Lin, J; Triki, F {\it Inverse scattering problems with multi-frequencies}. Inverse Problems {\bf 31} (2015), no. 9, 093001, 21 pp.

\bibitem{BLRX2015}
Bao, G; Lu, S; Rundell, W; Xu, B {\it A recursive algorithm for multifrequency acoustic inverse source problems}. SIAM J. Numer. Anal. {\bf 53} (2015), no. 3, 1608--1628.

\bibitem{BLT2010}
Bao, G; Lin, J; Triki, F {\it A multi-frequency inverse source problem}. J. Diff. Eq. {\bf 249} (2010), no. 12, 3443--3465.

\bibitem{BLZ2020}
Bao, G; Li, P; Zhao, Y {\it Stability for the inverse source problems in elastic and electromagnetic waves}. J. Math. Pures Appl. {\bf 134} (2020), 122--178.

\bibitem{BT2010}
Bao, G; Triki, F {\it Error estimates for the recursive linearization of inverse medium problems}. J. Comput. Math. {\bf 28} (2010), 725--744.

\bibitem{BT2020}
Bao, G; Triki, F {\it Stability for the multifrequency inverse medium problem}. J. Diff. Eq. {\bf 269} (2020), no. 9, 7106--7128.

\bibitem{C1980}
Calder\'{o}n, A P {\it On an inverse boundary value problem}. in Seminear on Numerical Analysis and Its Application to Continuum Physics, Rio de Jeneiro (1980), 65--73.


\bibitem{CFKKU2021}
C\^{a}rstea, C I; Feizmohammadi, A; Kian, Y; Krupchyk, K; Uhlmann, G {\it The Calder\'{o}n inverse problem for isotropic quasilinear conductivities}. Adv. Math. {\bf 391} (2021), 107956.

\bibitem{CIL2016}
Cheng, J; Isakov, V; Lu, S {\it Increasing stability in the inverse source problem with many frequencies}. J. Diff. Eq. {\bf 260} (2016), no. 5, 4786--4804.


\bibitem{EW2014}
Ev\'{e}quoz, G; Weth, T {\it Real Solutions to the Nonlinear Helmholtz Equation with Local Nonlinearity}. Arch. Ration. Mech. Anal. {\bf 211} (2014), no. 2, 359--388.


\bibitem{FO2020}
Feizmohammadi, A; Oksanen, L {\it An inverse problem for a semilinear elliptic equation in Riemannian geometries}. J. Diff. Eq. {\bf 269} (2020), 4683--4719.

\bibitem{FT2005}
Fibich, G; Tsynkov, S {\it Numerical solution of the nonlinear Helmholtz equation using nonorthogonal expansions}. J. Comput. Phys. {\bf 210} (2005), 183--224.


\bibitem{I2011}
Isakov, V {\it Increasing stability for the Schr\"{o}dinger potential from the Dirichlet-to-Neumann map}. Discrete Contin. Dyn. Syst. Ser. S {\bf 4} (2011), no. 3, 631--640.


\bibitem{IL2018}
Isakov, V; Lu, S {\it Increasing stability in the inverse source problem with attenuation and many frequencies}. SIAM J. Appl. Math. {\bf 78} (2018), no. 1, 1--18.

\bibitem{ILW2016}
\reviewC{Isakov, V; Lai, R-Y; Wang, J-N {\it Increasing stability for the conductivity and attenuation coefficients}. SIAM J. Math. Anal. {\bf 48} (2016), no. 1, 569--594.}

\bibitem{ILX2020}
Isakov, V; Lu, S; Xu, B {\it Linearized inverse Schr\"{o}dinger potential problem at a large wavenumber}. SIAM J. Appl. Math. {\bf 80} (2020), 338--358.

\bibitem{IW2014}
Isakov, V; Wang, J N {\it Increasing stability for determining the potential in the Schr\"{o}dinger equation with attenuation from the Dirichlet-to-Neumann map}. Inverse Probl. Imaging {\bf 8} (2014), no. 4, 1139--1150.

\bibitem{KKK2018}
\reviewA{Karamehmedovi\'{c}, M; Kirkeby, A; Knudsen, K {\it Stable source reconstruction from a finite number of measurements in the multi-frequency inverse source problem}. Inverse Problems {\bf 34} (2018), no. 6, 065004.}

\bibitem{KKU2020}
Kian, Y; Krupchyk, K; Uhlmann, G {\it Partial data inverse problems for quasilinear conductivity equations}. arXiv:2010.11409.

\bibitem{KLU2018}
Kurylev, Y; Lassas, M; Uhlmann, G {\it Inverse problems for Lorentzian manifolds and non-linear hyperbolic equations}. Invent. Math. {\bf 212} (2018), no. 3, 781--857.

\bibitem{KU2019}
Krupchyk, K; Uhlmann, G {\it A remark on partial data inverse problems for semilinear elliptic equations}. Proc. Amer. Math. Soc. {\bf 148} (2020), 681--685.

\bibitem{LLLS_I}
\reviewC{Lassas, M; Liimatainen, T; Lin, Y-H; Salo, M {\it Inverse problems for elliptic equations with power type nonlinearities}. J. Math. Pures Appl. {\bf 145} (2021), 44--82.}

\bibitem{LLLS_II}
Lassas, M; Liimatainen, T; Lin, Y-H; Salo, M {\it Partial data inverse problems and simultaneous recovery of boundary and coefficients for semilinear elliptic equations}. Rev. Mat. Iberoam. {\bf 37} (2021), no. 4, 1553--1580.

\bibitem{LLST2020}
\reviewA{Liimatainen, T; Lin, Y-H; Salo, M; Tyni, T {\it Inverse problems for elliptic equations with fractional power type nonlinearities}. J. Diff. Eq. {\bf 306} (2022), 189--219.}


\bibitem{NUW2013}
Nagayasu, S; Uhlmann, G; Wang, J {\it Increasing stability in an inverse problem for the acoustic equation}. Inverse Problems {\bf 29} (2013), no. 2, 025012, 11 pp.

\bibitem{WZ2018}
Wu, H; Zou, J {\it Finite element method and its analysis for a nonlinear Helmholtz equation with high wave numbers}. SIAM J. on Numer. Anal. {\bf 56} (2018), no. 3, 1338--1359.

\bibitem{XB2010}
Xu, Z; Bao, G {\it A numerical scheme for nonlinear Helmholtz equations with strong non-linear optical effects}. J. Opt. Soc. Am. A {\bf 27} (2010), 2347--2353.

\bibitem{YL2017}
Yuan, L; Lu, Y {\it Robust iterative method for nonlinear Helmholtz equation}. J. Comput. Phys. {\bf 343} (2017), 1--9.

\bibitem{ZZ2017}
\reviewA{Zhang, B; Zhang, H {\it Recovering scattering obstacles by multi-frequency phaseless far-field data}. J. Comput. Phys. {\bf 345} (2017), 58--73.}

\end{thebibliography}
\end{document}